\documentclass[journal,a4paper]{IEEEtran}

\usepackage[pdftex]{graphicx}
\DeclareGraphicsExtensions{.pdf,.jpeg,.png}

\usepackage[cmex10]{amsmath}
\usepackage{cite}
\usepackage{array}
\usepackage{fixltx2e}
\usepackage{url}

\usepackage{amsmath,amssymb}
\usepackage{amsfonts}

\DeclareMathOperator{\Ad}{Ad}
\DeclareMathOperator{\Tr}{Tr}

\DeclareMathOperator{\diag}{diag}
\DeclareMathOperator{\Aut}{Aut}

\begin{document}

\title{A Primer on Stochastic Differential Geometry for
Signal Processing}

\author{Jonathan~H.~Manton,~\IEEEmembership{Senior~Member,~IEEE}%
\thanks{Copyright (c) 2013 IEEE. Personal use of this material
is permitted. However, permission to use this material for any other
purposes must be obtained from the IEEE by sending a request to
\texttt{pubs-permissions@ieee.org}}%
\thanks{J. Manton is with the Control and Signal Processing Lab,
Department of Electrical and Electronic Engineering, The University of
Melbourne, Australia. e-mail: \texttt{j.manton@ieee.org}}%
\thanks{Manuscript received October 1, 2012; revised January 31, 2013
and March 31, 2013.}}

%

\maketitle

\newcommand{\ito}{It\^o}
\newcommand{\reals}{\mathbb{R}}
\newcommand{\Exp}{\operatorname{Exp}}
\newtheorem{proposition}{Proposition}
\newcommand{\cmplx}{\mathbb{C}}
\newcommand{\D}{\mathcal{D}}
\newcommand{\E}[1]{\mathbb{E}\left[{#1}\right]}
\newcommand{\Eg}[1]{\mathbb{E}^g\left[{#1}\right]}

\begin{abstract}
This primer explains how continuous-time stochastic processes
(precisely, Brownian motion and other \ito\ diffusions) can be
defined and studied on manifolds.  No knowledge is assumed of either
differential geometry or continuous-time processes.  The arguably
dry approach is avoided of first introducing differential geometry
and only then introducing stochastic processes; both areas are
motivated and developed jointly.  \end{abstract}

\begin{IEEEkeywords}
Differential geometry, 
stochastic differential equations on manifolds,
estimation theory on manifolds,
continuous-time stochastic processes,
\ito\ diffusions,
Brownian motion,
Lie groups.
\end{IEEEkeywords}

\IEEEpeerreviewmaketitle

\section{Introduction}
\IEEEPARstart{T}{he} tools of calculus --- differentiation,
integration, Taylor series, the chain rule and so forth --- have
extensions to curved surfaces and more abstract manifolds, and a
different set of extensions to stochastic processes.  Stochastic
differential geometry brings together these two extensions.

This primer was written from the perspective that, to be useful,
it should give more than a big-picture view by drilling down to
shed light on important concepts otherwise obfuscated in highly
technical language elsewhere.  Gaining intuition, and gaining the
ability to calculate, while going hand in hand, are distinct from
each other.  As there is ample material catering for the
latter~\cite{Clark:1973uv, Emery:1989cy, Rogers:2000ur, Protter:1986fj},
the focus is on the former.

Brownian motion plays an important role in both theoretical
and practical aspects of signal processing. Section~\ref{sec:two}
is devoted to understanding how Brownian motion can be defined on
a Riemannian manifold.  The standpoint is that it is infinitely
more useful to know how to simulate Brownian motion than to learn
that the generator of Brownian motion on a manifold is the
Laplace-Beltrami operator.

Stochastic development is introduced
early on, in Section~\ref{sec:bmsd}, because ``rolling
without slipping'' is a simple yet striking visual aid for understanding
curvature, the key feature making manifolds more complicated and
more interesting than Euclidean space.
Section~\ref{sec:ssm} explains how stochastic development can be
used to extend the concept of state-space models from Euclidean
space to manifolds.  This motivates the introduction of stochastic
differential equations in Section~\ref{sec:ito}.  Since textbooks
abound on stochastic differential equations in Euclidean space,
only a handful of pertinent facts are presented.  As explained in
Section~\ref{sec:sdem}, and despite appearances~\cite{Elworthy:1982tz,
Emery:1989cy, Hsu:2002tz}, going from stochastic
differential equations in Euclidean space to stochastic differential
equations on manifolds is relatively straightforward conceptually
if not technically.

Section~\ref{sec:acl} examines more closely the theory of stochastic
integration.  It explains (perhaps in a novel way) how randomness
can make it easier for integrals to exist.  It clarifies seemingly
contradictory statements in the literature about pathwise integration.
Finally, it emphasises that despite the internal complexities,
stochastic integrals are constructed in the same basic way as other
integrals and, from that perspective, are no more complicated than
any other linear operator.

The second half of the paper, starting with Section~\ref{sec:seom},
culminates in the generalisation of Gaussian random variables to
compact Lie groups and the re-derivation of the formulae
in~\cite{Said:2012hi} for estimating the parameters of these random
variables.  Particular attention is given to the special orthogonal
groups consisting of orthogonal matrices having unit determinant,
otherwise known as the rotation groups.  Symmetry makes Lie groups
particularly nice to work with.

Estimation theory on manifolds is touched on in Section~\ref{sec:afw},
the message being that an understanding of how an estimator will
be used is needed to avoid making \emph{ad hoc} choices about how
to assess the performance of an estimator, including what it means
for an estimator to be unbiased.

The reason for introducing continuous-time processes rather than
ostensibly simpler discrete-time processes is that the only linear
structure on manifolds is at the infinitesimal scale of tangent
spaces, allowing the theory of continuous-time processes to carry
over naturally to manifolds.  A strategy for working with discrete-time
processes is treating them as sampled versions of continuous-time
processes.  In the same vein, Section~\ref{sec:grv} uses Brownian
motion to generalise Gaussian random variables to Lie groups.

There are numerous omissions from this primer.  Even the \ito\
formula is not written down!  Perhaps the most regrettable is not
having the opportunity to explain why stochastic differentials
are genuine (second-order) differentials.

An endeavour to balance fluency and rigour has led to the label
\emph{Technicality} being given to paragraphs that may be skipped
by readers favouring fluency.  Similarly, regularity conditions
necessary to make statements true are routinely omitted.  \emph{Caveat
lector.}

\section{Simulating Brownian motion}
\label{sec:two}

\subsection{Background}
\label{sec:bm}

A continuous-time process in its most basic form is just an infinite
collection of real-valued random variables $X(t)$ indexed by time
$t \in [0,\infty)$.  Specifying directly a joint distribution for
an infinite number of random variables is generally not possible.
Instead, the following two-stage approach is usually adopted for
describing the statistical properties of a continuous-time process.

First and foremost, all the finite-dimensional joint distributions
of $X(t)$ are given; they determine most, but not all,
statistical properties of $X(t)$.  In detail, the finite-dimensional
joint distributions are the distributions of $X(\tau)$ for each
$\tau$, the pairwise joint distributions of $X(\tau_1)$ and $X(\tau_2)$
for all $\tau_1 \neq \tau_2$, and in general the joint distributions
of $X(\tau_1)$ to $X(\tau_n)$ for a finite but arbitrary $n$.

For fixed $\tau$, it is emphasised that $X(\tau)$ is simply a random
variable and should be treated as such; that $X(t)$ is a \textit{process}
is only relevant when looking at integrals or other limits involving
an infinite number of points.  However, the finite-dimensional
distributions on their own are inadequate for specifying the
distributions of such limits.  To exemplify, choose each $X(\tau)$
to be an \emph{independent} Gaussian random variable with zero mean
and unit variance, denoted $X(\tau) \sim N(0,1)$.  Although formally
a process, there is no relationship between any property of the
index set $[0,\infty)$ and any statistical property of the random
variables $X(t)$, $t \in [0,\infty)$.  For this process, $\lim_{t
\rightarrow \tau} X(t)$ and $\int_0^1 X(t)\,dt$ do not even
exist~\cite[p.\@45]{Wong:1985ky}.

Markov processes are examples of a relationship existing between
properties of the index set and statistical properties of the random
variables; a process is Markov if the distribution of any future
point $X(\tau+h)$, $h>0$, given past history $\{X(\tau_1),\cdots,X(\tau_n)
\mid \tau_1 < \cdots < \tau_n = \tau\}$, only depends on $X(\tau)$.
This memoryless property of Markov processes relates the ordering
of the index set $[0,\infty)$ to conditional independence of the
random variables.

Other examples are processes with continuous sample paths, where
the topology of the index set relates to convergence of random
variables.  In detail, if $X$ is a random variable then it is
customary to denote an outcome of $X$ by $x$.
Similarly, let $x(t)$ denote the realisation of a process
$X(t)$. When $x(t)$ is considered as the function $t \mapsto x(t)$,
it is called a sample path.  If (almost) all realisations $x(t)$
of a process $X(t)$ have continuous sample paths, meaning $t \mapsto
x(t)$ is continuous, then the process itself is called continuous.

The second step for defining a continuous-time process is describing
additional properties of the sample paths.  A typical example is
declaring that all sample paths are continuous.  Although the
finite-dimensional distributions do not define a process uniquely
--- so-called modifications are possible --- the additional requirement
that the sample paths are continuous ensures uniqueness.  (Existence
is a different matter; not all finite-dimensional distributions
are compatible with requiring continuity of sample paths.)

\noindent\textit{Technicality:}
For processes whose sample paths are continuous, the underlying
probability space~\cite{Williams:1991wk} can be taken to be
$(C([0,\infty)),\mathfrak{F},\mathcal{P})$ where $C([0,\infty))$
is the vector space of all real-valued continuous functions on the
interval $[0,\infty)$ and $\mathfrak{F}$ is the $\sigma$-algebra
generated by all cylinder sets. The probability measure $\mathcal{P}$
is uniquely determined by the finite-dimensional distributions of
the process.

If all finite-dimensional distributions are Gaussian then the process
itself is called a Gaussian process.  Linear systems
preserve Gaussianity; the output of a linear system driven by a
Gaussian process is itself a Gaussian process.  This leads to an
elegant and powerful theory of Gaussian processes in linear
systems theory, and is the theory often found in signal processing
textbooks.  Since manifolds are inherently nonlinear, such a
simplifying theory does not exist for processes on manifolds.
(Brownian motion can be defined on manifolds without reference to
Gaussian random variables.  Gaussian random fields can be defined
on manifolds, but these are real-valued processes indexed by a
manifold-valued parameter, as opposed to the manifold-valued processes
indexed by time that are the protagonists of this primer.)

\begin{figure}[!t]
\centering
\includegraphics[width=\columnwidth]{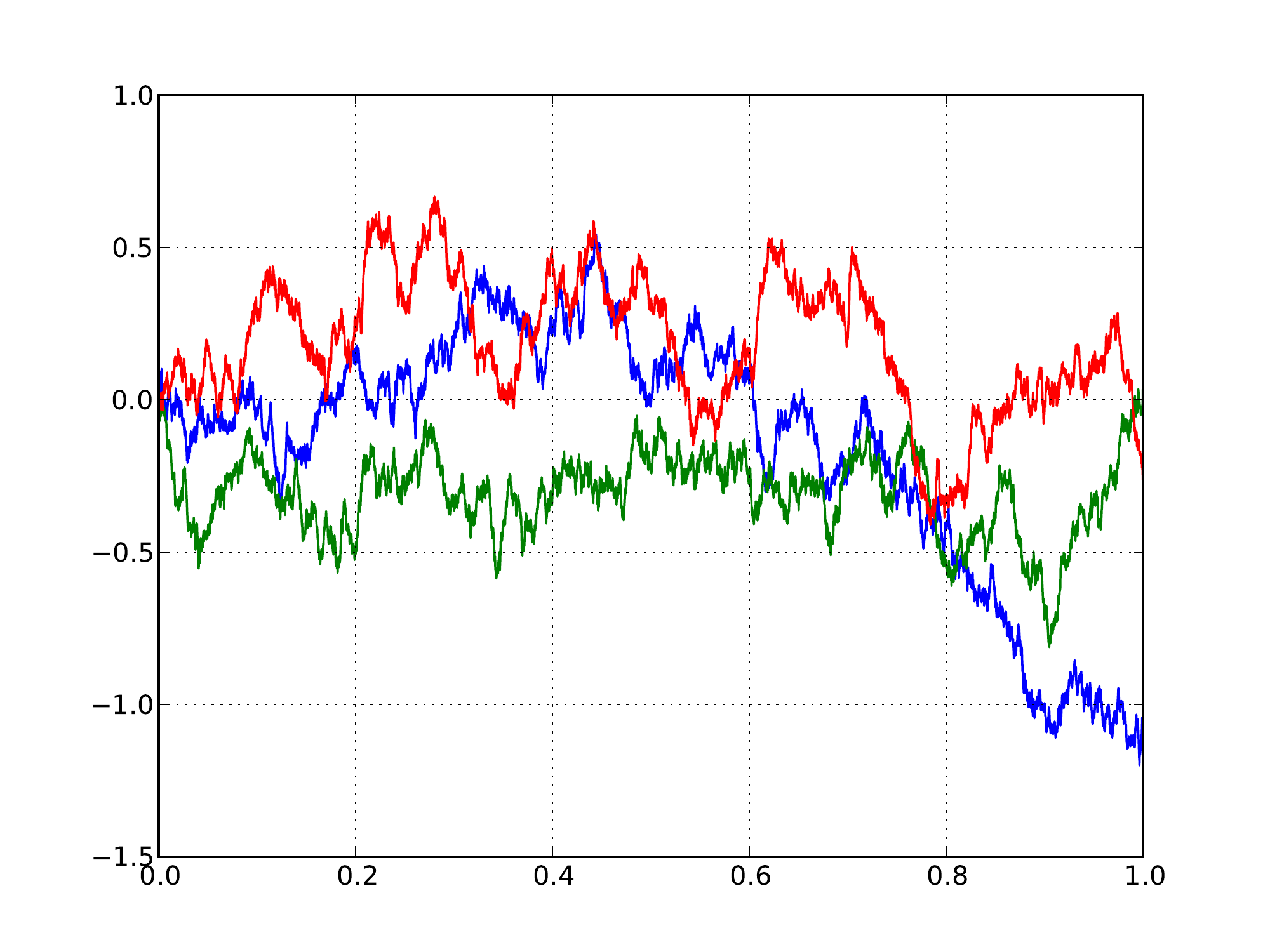}
\caption{Sample paths of Brownian motion on the closed
interval $[0,1]$ generated as described in the text using
$\delta t = 10^{-4}$.}
\label{fig:brown}
\end{figure}

The archetypal continuous Markov process is Brownian motion, normally
defined via its
finite-dimensional distributions and continuity of its sample
paths~\cite[Section 2.2]{bk:Oksendal:sde}.  In the
spirit of this primer though (and that of~\cite{Higham:2001fd}),
processes are best understood in the first instance by
knowing how to simulate them.  The sample paths of Brownian motion
plotted in Figure~\ref{fig:brown} were generated
as follows.  Set $X(0) = 0$ and fix a step size $\delta
t > 0$.  Note $\delta t$ is not the product of $\delta$ and $t$ but
merely the name of a positive real-valued quantity indicating a
suitably small discretisation of time.  Let $W(0),W(1),\cdots$ be
independent $N(0,1)$ Gaussian random variables.  Recursively define
$X((k+1)\,\delta t) = X(k\,\delta t) + \sqrt{\delta t}\,W(k)$ for
$k=0,1,2,\cdots$.  At non-integral values, define $X(t)$ by linear
interpolation of its neighbouring points: $X((k+\alpha)\,\delta
t) = (1-\alpha) X(k\,\delta t) + \alpha X((k+1)\,\delta t)$ for
$\alpha \in (0,1)$.  The generated process has the
correct distribution at integral sample points $t = 0, \delta t,
2\,\delta t, 3\,\delta t,\cdots$ and overall is an approximation
of Brownian motion converging (in distribution) to Brownian
motion as $\delta t \rightarrow 0$.

\noindent\textit{Technicalities:}
It can be shown that a process is
Brownian motion if and only if it is a \textit{continuous} process with
\textit{stationary} and \textit{independent} increments~\cite[p.\@
2]{Harrison:1990vt}; these properties force the process to be
Gaussian~\cite[Ch.\@ 12]{Breiman:1992vp}, a consequence of the
central limit theorem.  Brownian motion, suitably scaled and with
zero drift, is precisely the \emph{normalised}, or \emph{standard},
Brownian motion, described earlier.
Had the $W(k)$ been replaced by any other distribution with
zero mean and unit variance, the resulting process $X(t)$ still
would have converged to Brownian motion as $\delta t \rightarrow
0$.  Alternative methods for generating Brownian motion on the
interval $[0,1]$ include truncating the Karhunen-Lo\`eve expansion,
and successive refinements to the grid: $X(0)=0$, $X(1) \sim N(0,1)$,
and given neighbouring points $X(t)$ and $X(t+\delta
t)$, a mid-point is added by the rule $X(t + \frac{\delta t}2)
\sim \frac{X(t) + X(t+\delta t)}2 + \frac{\sqrt{\delta t}}2 N(0,1)$,
thus allowing $X(\frac12)$ to be computed with $\delta t = 1$, then
$X(\frac14)$ and $X(\frac34)$ with $\delta t = \frac12$ and so
forth.  Books specifically on Brownian motion include~\cite{
Revuz:1999wo, Chung:2002vc, Morters:2010vm}.
The origins of the mathematical concept of Brownian motion trace back
to three independent sources; Thiele (1880), Bachelier (1900) and
Einstein (1905).  According to~\cite{Jarrow:2004jh}, ``Of these
three models, those of Thiele and Bachelier had little impact for
a long time, while that of Einstein was immediately influential''.

\subsection{Brownian Motion and Stochastic Development}
\label{sec:bmsd}

Nothing is lost for the moment by treating manifolds as ``curved
surfaces such as the circle or sphere''.

Brownian motion models a particle bombarded randomly by much smaller
molecules.  The recursion $X((k+1)\,\delta t) = X(k\,\delta t) +
\sqrt{\delta t}\,W(k)$ introduced in Section~\ref{sec:bm} is thus
(loosely) interpreted as a particle being bombarded at regular time
instants.  Between bombardments, there is no force acting on the
particle, hence the particle's trajectory must be a curve of zero
acceleration. In Euclidean space, this implies particles move in
straight lines between bombardments, and explains why linear
interpolation was used earlier to connect $X(k\,\delta t)$ to
$X((k+1)\,\delta t)$.  On a manifold, a curve with zero acceleration
is called a geodesic.  Between bombardments, a particle on a manifold
travels along a geodesic.

Conceptually then, a piecewise approximation to Brownian motion on
a manifold can be generated essentially as before, just with
straight-line motion replaced by geodesic motion.

Since the Earth is approximately a sphere, long-distance travel
gives an intuitive understanding of the concepts of distance,
velocity and acceleration of a particle moving on the surface of a
sphere.  Travelling ``in a straight line'' on the Earth actually
means travelling along a great circle; great circles are the geodesics
of the sphere.

There are different ways of understanding geodesics, but the most
relevant for subsequent developments is the following: rolling a
sphere, without slipping, over a straight line drawn in wet ink on
a flat table will impart a curve on the sphere that precisely traces
out a geodesic.

\begin{figure}[!t]
\centering
\includegraphics[width=\columnwidth]{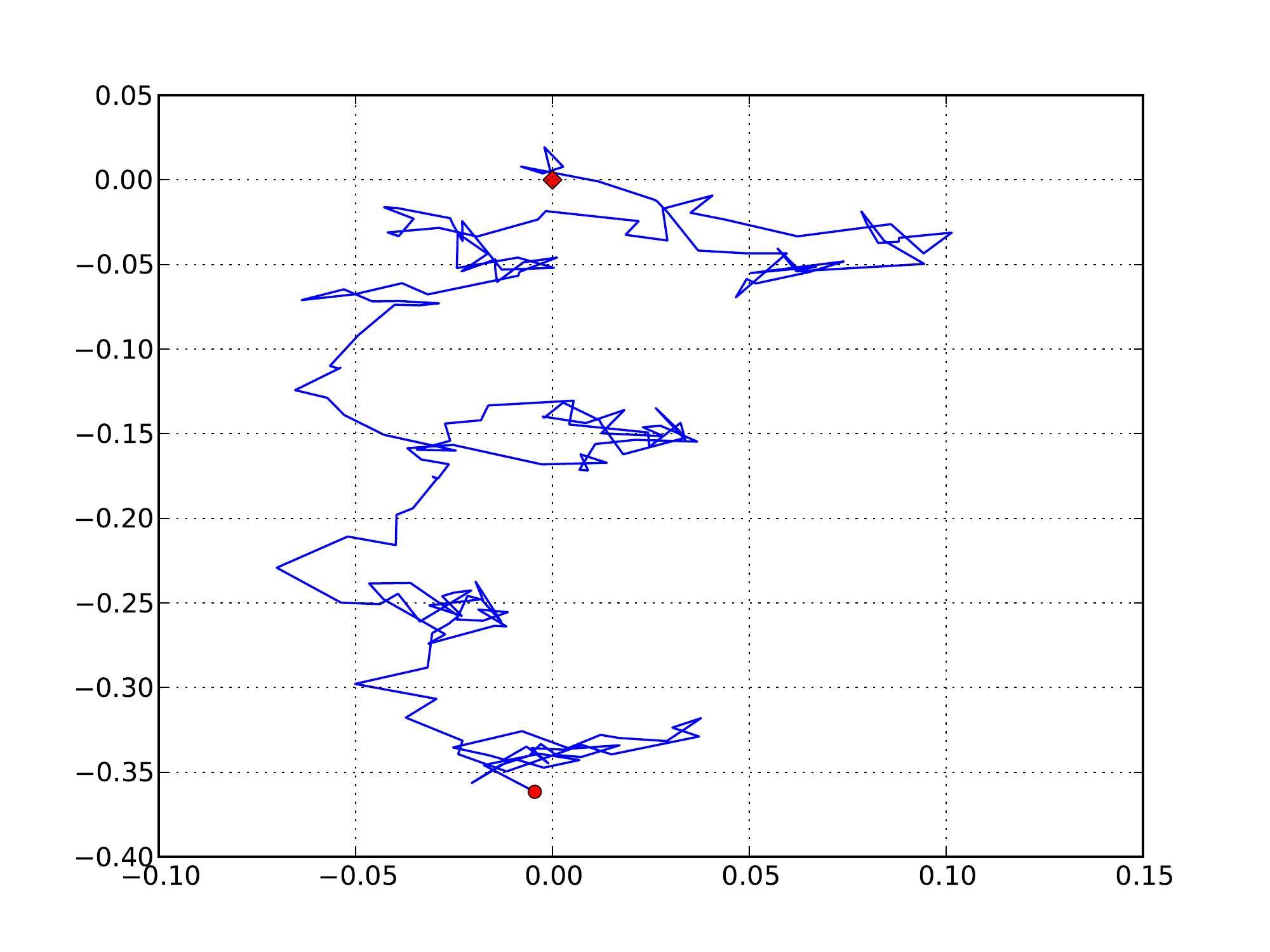}
\caption{Path traced out by a very short segment of a single
realisation of a two-dimensional Brownian motion starting from the
diamond at the origin and stopping at the circle; $X(t) =
(X_1(t),X_2(t))$ is two-dimensional Brownian motion if and only if
$X_1(t)$ and $X_2(t)$ are one-dimensional Brownian motions independent
of each other.  Time has been omitted.}
\label{fig:brown2D}
\end{figure}

Rolling a manifold over a piecewise smooth curve in Euclidean space
to obtain a curve on the manifold is called \textit{development}.
The development of a piecewise linear curve is a piecewise geodesic
curve.  If $x(t) = (x_1(t),x_2(t))$ is the piecewise linear
approximation in Figure~\ref{fig:brown2D} of a realisation of
two-dimensional Brownian motion then the developed piecewise geodesic
curve obtained by rolling the sphere over the curve in
Figure~\ref{fig:brown2D} is a piecewise geodesic approximation of
a realisation of Brownian motion on the sphere.

Development is defined mathematically as the solution of a certain
differential equation.  If
a curve is not differentiable then it cannot be developed in the classical
sense. A sample path of Brownian motion is nowhere differentiable!
Therefore, in the first instance, development can only be used to
take a piecewise smooth \textit{approximation} of Brownian motion
in Euclidean space and obtain a piecewise smooth \textit{approximation}
of Brownian motion on a manifold.  Nevertheless, it is possible to
``take limits'' and develop a theory of \textit{stochastic development}.
The stochastic development of Brownian motion in Euclidean space
is the limiting process obtained by developing successively more
accurate piecewise smooth approximations of the Brownian motion in
Euclidean space.

\noindent\textit{Technicality:}
This ``smooth approximation'' approach to stochastic development
is more tedious to make rigorous than the stochastic differential
equation approach taken in~\cite[Chapter 2]{Hsu:2002tz} but offers
more intuition for the neophyte.

The inverse of development is anti-development and can be visualised
as drawing a curve on a manifold in wet ink and rolling the manifold
along this curve on a table thereby leaving behind a curve
on the table.  Its stochastic counterpart, stochastic anti-development,
can be thought of as the limiting behaviour of the anti-development
of piecewise geodesic approximations of sample paths.

As a rule of thumb, techniques (such as filtering~\cite{Wonham:1965wk})
for processes on Euclidean space can be extended to processes on
manifolds by using stochastic anti-development to convert the
processes from manifold-valued to Euclidean-space-valued. Although
this may not always be computationally attractive, it nevertheless
affords a relatively simple viewpoint.

A process on a manifold is Brownian motion if and only if its
stochastic anti-development is Brownian motion in Euclidean
space~\cite[(8.26)]{Emery:1989cy}.  
To put this in perspective, note Brownian motion on a manifold
cannot be defined using finite-dimensional distributions because
there is no \emph{direct} definition of a Gaussian random variable
on a manifold.  (Even disregarding the global topology of a manifold,
a random variable which is Gaussian with respect to one chart need
not be Gaussian with respect to another chart.) Much of the
Euclidean-space theory relies on linearity, and the only linearity
on manifolds is at the infinitesimal scale of tangent
spaces.  Stochastic development operates at the infinitesimal scale,
replicating as faithfully as possible on a manifold a process in
Euclidean space.

Although Gaussian random variables cannot be used to define Brownian
motion on a manifold, the reverse is possible; Brownian motion can
be used to generalise the definition of a Gaussian random variable
to a manifold; see Section~\ref{sec:grv}.

\noindent\textit{Technicalities:}
There are a number of characterisations of Brownian motion that can
be used to generalise Brownian motion to manifolds, including as
the unique \ito\ diffusion generated by the Laplace operator.  The
end result is nevertheless the same~\cite{Meyer:1981vy, Elworthy:1982tz,
Lewis:1986hi, Emery:1989cy, Ikeda:1989wt, Rogers:2000ur, Hsu:2002tz}.
Whereas for deterministic signal processing, and in particular for
optimisation on manifolds, there is benefit in not endowing a
manifold with a Riemannian metric~\cite{Manton:opt_mfold, Manton:2012vk},
Brownian motion must be defined with respect to a Riemannian metric.
Although some concepts, such as that of a semimartingale, can be
defined on a non-Riemannian manifold, it is simplest here to assume
throughout that all manifolds are \textbf{Riemannian manifolds}.

\subsection{The Geometry of the Sphere}
\label{sec:sphere}

This section continues the parallel threads of stochastic development
and Brownian motion.  Equations are derived for simulating Brownian
motion on a sphere by rolling a sphere along a simulated path of
Brownian motion in $\reals^2$.  It is convenient to change reference
frames and roll a sheet of paper around a sphere than roll a sphere
on a sheet of paper.

\noindent\textit{Exercise:}
Mentally or otherwise, take a sheet of graph paper and a soccer
ball.  Mark a point roughly in the middle of the graph paper as
being the origin, and draw two unit-length vectors at right-angles
to each other based at the origin.  Place the origin of the graph
paper on top of the ball.  Roll the paper down the ball until
arriving at the equator.  Then roll the paper along the equator for
some distance, then roll it back up to the top of the ball.  Compare
the orientation of the two vectors at the start and at the end of
this exercise; in general, the orientation will have changed due
to the curvature of the sphere.  (Furthermore, in general, the
origin of the graph paper will no longer be touching the top of the
sphere.)

Throughout, any norm $\|\cdot\|$ or inner product $\langle \cdot,\cdot
\rangle$ on $\reals^n$ is the Euclidean norm or inner product.
Perpendicular vectors are denoted $v \perp p$, meaning $\langle v,p
\rangle = 0$.

Take the sphere $S^2 = \{(x_1,x_2,x_3) \in \reals^3
\mid x_1^2+x_2^2+x_3^2=1\}$.  Start the Brownian motion at the North
pole: $B(0) = (0,0,1)$.  Place a sheet of graph paper on top of the
sphere, so the origin $(0,0)$ of the graph paper makes contact with
the North pole $(0,0,1)$.  Let $w_1(0)$ and $w_2(0)$ be realisations
of independent $N(0,1)$ random variables.  On the graph paper, draw
a line segment from the origin $(0,0)$ to the point $(\sqrt{\delta
t}\,w_1(0), \sqrt{\delta t}\, w_2(0))$; recall from Section~\ref{sec:bm}
and Figure~\ref{fig:brown2D} that this is the first segment of an
approximate sample path of Brownian motion on $\reals^2$.  The paper
is sitting inside $\reals^3$, and the point $(\sqrt{\delta t}\,w_1(0),
\sqrt{\delta t}\,w_2(0))$ on the paper is actually located at the
point $(\sqrt{\delta t}\,w_1(0), \sqrt{\delta t}\,w_2(0),1)$ in
$\reals^3$ because the paper is lying flat on top of the sphere.

Rolling the paper down the sphere at constant velocity along the
line segment so it reaches the end of the segment at time $\delta t$
results in the point of contact between the paper and
the sphere being given by
\begin{equation}
B(t) = \cos\left(\frac{t}{\delta t} \|v\|\right) p +
    \sin\left(\frac{t}{\delta t} \|v\|\right) \frac{v}{\| v\|}
\end{equation}
for $t \in [0,\delta t]$, where $p=(0,0,1)$ is the original point
of contact (the North pole) and $v = (\sqrt{\delta t}\,w_1(0),
\sqrt{\delta t}\,w_2(0),0)$.
This can be derived from first principles by requiring $B(t)$ to
remain in a plane and have constant angular velocity.

In general, given any
$p \in S^2$ and $v \perp p$, let $\Exp_p(v)$ denote the final
contact point of the sphere and piece of paper obtained by
starting with the paper touching the sphere at $p$, marking
on the paper a line segment from $p$ to $p+v$, and rolling
the paper over the sphere along that line.  (Since $v \perp p$,
and the paper is tangent to the sphere at $p$,
the point $p+v$ will lie on the paper.) The curve
$t \mapsto \Exp_p(tv)$ is a geodesic and follows a great circle.
Explicitly, $\Exp_p(0) = p$ and, for $\|v\| \neq 0$,
\begin{equation}
\Exp_p(v) = \cos\left(\|v\|\right) p +
    \sin\left( \|v\|\right) \frac{v}{\| v\|}.
\end{equation}
At the risk of belabouring the point, if $S^2$ represents the Earth
and a person sets out from a point $p$ with initial velocity $v$
and continues ``with constant velocity in a straight line'' then
his position at time $t$ will be $\Exp_p(tv)$.  In detail, each
step actually involves first taking a perfectly straight step in
the tangent plane, meaning the leading foot will be slightly off the
Earth, then without swivelling on the other foot, letting the force
of gravity pull the leading foot down to the closest point on Earth.
This is a discrete version of ``rolling without slipping'' and hence
produces (or defines) a geodesic in the limit as smaller and smaller
steps are taken.  

The Riemannian exponential map $\Exp_p(v)$ can be defined analogously
on any Riemannian manifold.  The set of allowable velocities $v$
for which $\Exp_p(v)$ makes sense is called the tangent space to
the manifold at the point $p$; just like for the sphere, the tangent
space can be visualised as a sheet of paper providing the best
linear approximation to the shape of the manifold in a neighbourhood
of the point of contact and $v$ must be such that $p+v$ lies on
this (infinitely large) sheet of paper.  Alternatively, if a
sufficiently small (or sufficiently short-sighted) ant were standing
on the manifold at point $p$, so that the manifold looked flat,
then the set of possible directions (with arbitrary magnitudes) the
ant could set out in from his perspective forms the tangent space
at $p$.

\noindent\textit{Technicalities:}
The Riemannian exponential function is defined via a differential
equation.  If the manifold is not complete, the differential equation
may ``blow up''; this technicality is ignored throughout the primer.
Since every Riemannian manifold can be embedded in a sufficiently
high-dimensional Euclidean space, this primer assumes for simplicity
that all Riemannian manifolds are subsets of Euclidean space. The
Riemannian geometry of such a manifold is determined by the Euclidean
geometry of the ambient space; the Euclidean inner product induces
an inner product on each tangent space.  This is consistent with
defining the length of a curve on a manifold as the Euclidean
length of that curve when viewed as a curve in the ambient Euclidean
space.

Returning to simulating Brownian motion on the sphere, recall the
original strategy was to generate a Brownian motion on the plane
then develop it onto the sphere.  Carrying this out exactly would
involve keeping track of the orientation of the paper as it moved
over the sphere.  Although this is easily done, it is not necessary
for simulating Brownian motion because Gaussian random vectors are
symmetric and hence invariant with respect to changes in orientation.

If a particle undergoing Brownian motion is currently
at the point $p \in S^2$, its next position, after $\delta t$
units of time, can be simulated by generating a three-dimensional
Gaussian random vector $v \sim N(0,I) \in \reals^3$, projecting $v$
onto $T_pS^2$ --- replace $v$ by $v - \langle v,p \rangle p$ ---
and declaring the next position of the particle to be $\Exp_p(\sqrt{\delta
t}\, v)$.  This generalises immediately to arbitrary manifolds and
is summarised in Section~\ref{sec:Bmm}.  (Alternatively, given an
orthonormal basis for $T_pS^2$, a two-dimensional Gaussian random
vector could have been used to generate an appropriate random element
of $T_pS^2$.)

Although orientation was ultimately not needed for defining Brownian
motion, it is an important concept by which to understand curvature,
and enters the picture for more general processes (such as Brownian
motion with drift).

\noindent\textit{Technicality:}
For a non-embedded manifold $M$, the natural setting for (stochastic)
development is the frame bundle $\mathcal{F}(M)$ of $M$ equipped
with a connection~\cite{Hsu:2002tz, Elworthy:1982tz, Rogers:2000ur}.
The connection decomposes the tangent bundle of
$\mathcal{F}(M)$ into horizontal and
vertical components, and leads to the concept of a horizontal
process.  A horizontal process on $\mathcal{F}(M)$ is essentially
a process on $M$ augmented by its current orientation.  If
$M$ is Riemannian then the orthonormal frame bundle $\mathcal{O}(M)$
can be used instead of $\mathcal{F}(M)$.  Horizontal Brownian motion
can be defined on $\mathcal{O}(M)$ via a Stratonovich stochastic
differential equation that stochastically develops Brownian motion
in Euclidean space onto the horizontal component (with respect to
the Levi-Civita connection) of $\mathcal{O}(M)$. The bundle projection
of this horizontal Brownian motion yields Brownian motion on $M$.

\subsection{A Working Definition of a Riemannian Manifold}
\label{sec:rm}

For the purposes of this primer, manifolds are defined as subsets
of Euclidean space that are sufficiently nice to permit a useful
theory of differentiation of functions from one such subset to
another.
(Furthermore, only $C^\infty$-smooth manifolds are discussed.)
Conditions will be given for a subset $M \subset \reals^n$
to be an $m$-dimensional manifold for some positive integer $m \leq
n$.  (A zero-dimensional manifold is a countable collection of
isolated points and will not be considered further.) 
This will confirm the circle $S^1$ and
sphere $S^2$ as manifolds of dimension one and two, respectively.
Graphs of smooth functions are prototypical manifolds:
if $f\colon \reals^m \rightarrow \reals^{n-m}$ is a smooth function,
meaning derivatives of all orders exist, its graph $M =
\{(x,f(x)) \in \reals^m \times \reals^{n-m} \cong \reals^n \mid x
\in \reals^m\}$ is an $m$-dimensional manifold. 

For each point $p \in
M$, define $T_pM$ as the set of all possible velocity vectors
$\gamma'(0)$ taken on by smooth 
curves $\gamma\colon \reals \rightarrow
\reals^n$ whose images are wholly contained in $M$ and that pass
through $p$ at time $0$.  For example,
if $p \in M = S^2$ then the (only) requirements on $\gamma$ are
that it is smooth, that $\gamma(0)=p$ and $\|\gamma(t)\|=1$.
In symbols,
\begin{equation}
T_pM = \{\gamma'(0) \mid \gamma\colon \reals \rightarrow \reals^n,
    \, \gamma(0)=p,\, \gamma(\reals) \subset M\}
\end{equation}
where it is implicitly understood that $\gamma$ must be infinitely
differentiable.  (No difference results if $\gamma$ is only
defined on an open neighbourhood of the origin; usually such
curves are denoted $\gamma\colon (-\epsilon,\epsilon) \rightarrow M$
in the literature.)

The first requirement placed on $M$ is for it to look infinitesimally 
like $\reals^m$.  Precisely, it is required that $T_pM \subset
\reals^n$ is an $m$-dimensional vector subspace of $\reals^n$ for
every $p \in M$.  This prevents $M$ from having (non-tangential)
self-intersections, e.g., the letter $X$ is not a manifold, and it
prevents $M$ from having cusps, e.g., the letter $V$ is not a
manifold because no \textit{smooth} curve $\gamma$ passes through
the bottom tip of $V$ except for the constant curve with $\gamma'(0)=0$.

Usually this first requirement is not stated because it
is subsumed by requiring the manifold be locally Euclidean,
defined presently.  Nevertheless, it emphasises the importance of
tangent spaces.  The visual image of a piece of paper placed against
a sphere at the point $p \in S^2 \subset \reals^3$ is the \emph{affine}
tangent space.  The tangent space $T_pS^2$ is obtained by taking
the piece of paper and translating it to pass through the origin
of the ambient space $\reals^3$.  This distinction is usually
blurred.

The first requirement fails to disqualify the figure of eight
from being a manifold because it cannot detect
\textit{tangential} self-intersections.  This can only be detected
by considering what is happening in a neighbourhood of each and
every point; 
it is required that for all $p \in M \subset \reals^n$
there exists a sufficiently small open ball
$B_r(p) = \{z \in \reals^n \mid \|z-p\| < r\}$ of
radius $r > 0$ and a diffeomorphism
$h\colon \reals^n \rightarrow B_r(p)$, meaning $h$ is bijective and
both $h$ and its inverse are smooth, such that
\begin{equation}
\label{eq:h}
h(\reals^m \times \{0\}) = B_r(p) \cap M
\end{equation}
where $\reals^m \times \{0\} = \{(x,0) \in \reals^n \mid x \in \reals^m\}$
is the embedding of $\reals^m$ into $\reals^n$ obtained by setting
the last $n-m$ coordinates equal to zero.
(The basic intuition is that the classical calculus on a flat subspace
$\reals^m \times \{0\}$ of $\reals^n$ should be extendable to a
calculus on diffeomorphic images of this flat subspace.)

The restriction of $h$ in (\ref{eq:h}) to $\reals^m \times \{0\}$
is a parametrisation of a part of the manifold $M$, however, the
condition is stronger than this since it requires the parametrisation
include \textit{all} points of $M$ in $B_r(p) \cap M$ and \textit{no
other} points.  This excludes the figure ``8'' because at the point
where the top and bottom circles meet, every one-dimensional
parametrisation can get at best only \textit{part} of the lower
hemisphere of the top circle and the upper hemisphere of the bottom
circle.  In fact, (\ref{eq:h}) ensures that every manifold locally
looks like a rotated graph of a smooth function.  

A manifold $M \subset \reals^n$ inherits a Riemannian structure
from the Euclidean inner product on $\reals^n$.  This leads to
defining the acceleration of a curve $\gamma\colon \reals \rightarrow M
\subset \reals^n$ at time $t$ as
$\pi_{\gamma(t)}\left(\gamma''(t)\right)$ where $\pi_p\colon \reals^n
\rightarrow T_pM$ is orthogonal projection onto $T_pM$.  The curve
$\gamma$ is a geodesic if and only if $\gamma''(t)$ contains only
those vectorial components necessary to keep the curve on the
manifold, that is, $\pi_{\gamma(t)}\left(\gamma''(t)\right) = 0$
for all $t$.  This same condition could have been deduced 
by developing a straight line onto $M$, as in Section~\ref{sec:bmsd}. 

\noindent\textit{Technicalities:}
On a non-Riemannian manifold, there is \textit{a priori} no way of
defining the acceleration of a curve because there is no
\textit{distinguished} way of aligning $T_pM$ and $T_qM$ for two
distinct points $p$ and $q$.  The above definition implicitly uses
the Riemannian structure coming from the ambient space $\reals^n$
and accords with comparing tangent vectors at two distinct points
of a curve by placing a piece of paper over the manifold at point
$q$, drawing the tangent vector at $q$ on the piece of paper, then
rolling the paper to $p$ along the curve, and drawing the tangent
vector at $p$ on the paper.  Because the piece of paper represents
Euclidean space, the base points of the vectors drawn on the paper
can be translated in the usual way so that they align.  This then
allows the difference of the two vectors to be taken, and ultimately,
allows the acceleration to be defined as the rate of change of the
velocity vectors.  A mechanism for aligning neighbouring affine
tangent spaces along a curve is called an \textit{affine connection}.
The particular affine connection described here is the Levi-Civita
connection, the unique torsion-free connection that is compatible
with the Riemannian metric.  In more advanced settings, there may
be advantages to using other affine connections.  A limitation of
insisting manifolds are subsets of Euclidean space is that changing
to a different metric requires changing the embedding,
for example, changing the sphere into an ellipse.

\subsection{Brownian Motion on Manifolds}
\label{sec:Bmm}

Assembling the pieces leads to the following algorithm for simulating
Brownian motion on an $m$-dimensional Riemannian manifold $M \subset
\reals^n$.

Choose a starting point on $M$; set $B(0)$ to this point.  Fix a
step size $\delta t > 0$.  For $k=0,1,\cdots$, recursively do the
following.  Generate a Gaussian random vector $W(k) \in T_{B(k\,\delta
t)}M \subset \reals^n$, either with the help of an orthonormal basis
for $T_{B(k\,\delta t)}M$, or by generating an $n$-dimensional
$N(0,I)$ Gaussian random vector and projecting the vector orthogonally
onto $T_{B(k\,\delta t)}M$ to obtain $W(k)$.
Then define
\begin{equation}
\label{eq:B}
B(k\,\delta t + t) = \Exp_{B(k\,\delta t)}\left( \frac{t}{\delta t}
    \sqrt{\delta t}\,W(k) \right).
\end{equation}
for $t \in [0,\delta t]$.
If $M = \reals^n$ then $\Exp_p(v) = p+v$ and
(\ref{eq:B}) agrees with the algorithm in Section~\ref{sec:bm}.

\noindent\textit{Technicality:}
An advantage of projecting orthogonally onto the tangent space
rather than constructing an orthonormal basis is that while the
orthogonal projection $\pi_p\colon \reals^n \rightarrow T_pM$ varies
smoothly in $p$, on non-parallelisable manifolds it is not possible
to find a \emph{continuous} mapping from $p \in M$ to an orthonormal
basis of $T_pM$.  The hairy ball theorem implies that the sphere
$S^2$ is not parallelisable.  (In fact, in terms of spheres, only
$S^0$, $S^1$, $S^3$ and $S^7$ are parallelisable.)

It is verified
in~\cite{Gangolli:1964gm} that, as $\delta t \rightarrow 0$, the
above approximation converges in distribution to Brownian motion,
where Brownian motion is defined by some other means.
(For the special case of Lie groups, see also~\cite{McKeanJr:1960vj}.
For hypersurfaces, see~\cite{vandenBerg:1985jj}.)
Nevertheless, engineers (and physicists) may find it attractive to
treat the limit of (\ref{eq:B}) as the definition of Brownian motion.
All (\ref{eq:B}) is saying is that at each step, movement in any
direction is equally likely and independent of 
the past.  By the central limit theorem, it suffices for the
$W(k)$ to have zero mean and unit variance; see~\cite{Binder:2009vk}
for an analysis of an algorithm commonly used in practice. The
presence of the square root in the term $\sqrt{\delta t}$ is easily
explained by the compatibility requirement that the variance of
$B(T)$ generated using a step size of $\delta t$ be equal to the
variance of $B(T)$ generated using a step size of $\frac{\delta
t}2$; if this were not so then the processes need not converge as
$\delta t \rightarrow 0$.

While (\ref{eq:B}) is suitable for numerical work, for calculations
by hand it is convenient to ``take limits'' and work with the actual
process.  A direct analogy is preferring to work with $\frac{dx}{dt}$
rather than $\frac{x(t+\delta t) - x(t)}{\delta t}$.  Section~\ref{sec:ito}
introduces a stochastic calculus.

\section{State-Space Models on Manifolds}
\label{sec:ssm}

\subsection{Motivation}
\label{sec:ssmmot}

Signal processing involves generating new processes from old.  In
Euclidean space, a process can be passed through a linear
time-invariant system to obtain a new process.  This can be written
in terms of an integral and motivates asking if a continuous-time
process evolving on a manifold, such as Brownian motion, can be
integrated to obtain a new process.

Another obvious question is how to generalise to manifolds state-space
models with additive noise.  The classical linear discrete-time
state-space model is
\begin{align}
\label{eq:sse}
X_{k+1} &= A_k X_k + B_k V_k, \\
\label{eq:oe}
Y_k &= C_k X_k + D_k W_k
\end{align}
where $A_k, B_k, C_k, D_k$ are matrices and $V_k,W_k$ are random
vectors (noise).  The vector $X_k$ is the state at time $k$,
and the state-space equation represents the dynamics governing how
the state changes over time.  It comprises a deterministic part
$A_k X_k$ and a stochastic part $B_k V_k$.  The second equation is
the observation equation: the only measurement of the state available
at time $k$ is the vector $Y_k$ comprising a linear function 
of the state and additive noise.

There does not appear to be a natural generalisation of
\emph{discrete-time} state-space models to arbitrary manifolds
because it is not clear how to handle the addition of the two terms
in each equation.  (In some cases, group actions could be used.)
It will be seen presently that \emph{continuous-time} state-space
models generalise more easily.  This suggests the expediency of
treating discrete-time processes as sampled versions of continuous-time
processes.

The continuous-time version of (\ref{eq:sse}) would be
\begin{equation}
\label{eq:ctsse}
\frac{dX}{dt} = A(t) X(t) + B(t) \frac{dV(t)}{dt}
\end{equation}
if the noise process $V(t)$ was sufficiently nice that its sample
paths were absolutely continuous; that this generally is not 
the case is ignored for the moment.

Although $Y(t) = C(t) X(t) + D(t) W(t)$ is an analogue of (\ref{eq:oe}),
usually the observation process takes instead the form
\begin{equation}
\label{eq:ctoe}
\frac{dY}{dt} = C(t) X(t) + D(t) \frac{dW(t)}{dt}
\end{equation}
which integrates rather than instantaneously samples the state.

Although the right-hand sides of (\ref{eq:sse}) and (\ref{eq:ctsse})
are sums of two terms, crucially, it is two tangent vectors being
summed in (\ref{eq:ctsse}).  Two points on a manifold cannot be
added but two tangent vectors in the same tangent space can.
Therefore, (\ref{eq:ctsse}) and (\ref{eq:ctoe}) extend
naturally to manifolds provided the terms $A(t)X(t)$ and $C(t)X(t)$
are generalised to be of the form $b(t,X(t))$; see (\ref{eq:Xb}).
The challenge is if $V(t)$ and $W(t)$ are
Brownian motion then (\ref{eq:ctsse}) and (\ref{eq:ctoe}) require
an appropriate interpretation because Brownian motion is nowhere
differentiable (almost surely).  The following subsections hint
at how this is done via piecewise approximations, and Section~\ref{sec:ito}
gives a rigorous interpretation by changing (\ref{eq:ctsse}) and
(\ref{eq:ctoe}) to integral form.

\subsection{Modelling the State Process}

Building on the material in Sections~\ref{sec:two} and~\ref{sec:ssmmot},
an attempt is made to model a particle moving on a sphere subject
to disturbance by Brownian motion.  Let $X(t) \in S^2 \subset
\reals^3$ denote the position of the particle at time $t$.  Its
deterministic component can be specified by a differential equation
\begin{equation}
\label{eq:Xb}
\frac{dX}{dt} = b(t,X(t)).
\end{equation}
Provided $b(t,X(t))$ lies in the tangent space of the sphere at $X(t)$,
the solution of (\ref{eq:Xb}) is forced to lie on the sphere if
the initial point $X(0)$ does.
A simple numerical solution of (\ref{eq:Xb}) is obtained by
combining the forward-Euler method with the Riemannian exponential
function, the latter ensuring the approximate solution 
remains on the sphere:
\begin{equation}
X(t+\delta t) = \Exp_{X(t)}\big( \delta t\,b(t,X(t)) \big).
\end{equation}
Referring to (\ref{eq:B}) with $t=\delta t$, the following idea
presents itself:
\begin{equation}
\label{eq:Xsimple}
X(t+\delta t) = \Exp_{X(t)}\big( \delta t\,b(t,X(t))
    + \sqrt{\delta t}\,W(t)\big)
\end{equation}
where $W(t)$ is an $N(0,I)$ Gaussian random vector in $\reals^3$
projected orthogonally onto the tangent space $T_{X(t)}S^2$.  
The (approximately instantaneous) velocity of the particle
at time $t$ is the sum of a deterministic component and a random
component.  Continuous-time approximations can be obtained by
interpolating using geodesics, as in (\ref{eq:B}), in which case
(\ref{eq:Xsimple}) converges to a well-defined process
on the sphere~\cite{Duncan:1976vx}.

\subsection{Modelling the Observation Process}

Notwithstanding that the first two cases are subsumed by the third,
the three cases of interest are: the state process evolves
in Euclidean space yet the observation process is manifold-valued;
the state process evolves on a manifold but the observation process
is real-valued; and, the state and observation processes are
manifold-valued.

If the state process $X(t)$ evolves in Euclidean space then
stochastic development can be used to feed it into the
observation process~\cite{Duncan:1976vx}:
\begin{equation}
\label{eq:YobsI}
Y(t+\delta t) = \Exp_{Y(t)}\big( X(t+\delta t) - X(t)
    + \sqrt{\delta t}\,W(t)\big).
\end{equation}
If $X(t)$ is not observed directly, but only $g(X(t))$
where $g$ is a
smooth function between Euclidean spaces, then
a straightforward modification of (\ref{eq:YobsI}) is
\begin{equation}
Y(t+\delta t) = \Exp_{Y(t)}\big( g(X(t+\delta t)) - g(X(t))
    + \sqrt{\delta t}\,W(t)\big).
\end{equation}
In other words, first the new process $\tilde X(t) = g(X(t))$ is
formed, then noise is added to it, and finally it is stochastically
developed (via piecewise approximations) onto the manifold.

If the state process $X(t)$ evolves on a manifold of
dimension $m$ but the observation
process $Y(t)$ is real-valued then stochastic anti-development can
be used~\cite{Darling:1984ca}.  In a sufficiently small domain,
$\Exp$ is invertible, and the instantaneous velocity of $X(t)$,
which in general does not exist, can nevertheless be approximated by
$\Exp^{-1}_{X(t)}\big(X(t+\delta t) - X(t)\big)$.
This produces a vector in $\reals^m$ which can be used
to update an observation process evolving in Euclidean space.
(By interpreting differential equations as integral equations, as
discussed in Section~\ref{sec:ito}, neither
$X(t)$ nor $Y(t)$ need be differentiable for there to be a
well-defined limiting relationship between the instantaneous velocities of
piecewise approximations of the processes.)

Finally, the general case of both $X(t)$ and $Y(t)$ evolving on
manifolds falls under the framework of stochastic differential
equations between manifolds~\cite[Section 3]{Cohen:2000hp}.  
Basically,
$\Exp^{-1}_{X(t)}$ can be used to obtain a real-valued vector
approximating the instantaneous velocity of $X(t)$, which after a
possible transformation, can be fed into $\Exp_{Y(t)}$ to update
the observation process $Y(t)$.  See~\cite{Emery:1989cy} for
details.

\subsection{Discussion}

Section~\ref{sec:ssm} demonstrated, at least on an intuitive
level, that stochastic development and anti-development, and the
$\Exp$ map in particular, provide a relatively straightforward way
of generalising state-space models to manifolds.

Nevertheless, it is important to understand what the limiting
processes are that $\Exp$ is being used to approximate.  This is
the purpose of stochastic calculus, and once it is appreciated that
the ``smooth approximation'' approach discussed in this primer leads
to a stochastic calculus, it is generally easier to use directly
the stochastic calculus.

\section{Stochastic Calculus and \ito\ Diffusions}
\label{sec:ito}

This section does not consider manifolds or 
processes with jumps~\cite{Applebaum:2009ix}.  
Standard references include~\cite{Durrett:1996wh, bk:Oksendal:sde}.

\subsection{Background}
\label{sec:detan}

A generalisation of (\ref{eq:Xb}) is the \textit{functional equation}
\begin{equation}
\label{eq:Xbint}
X(t) = X(0) + \int_0^t b(s,X(s))\,ds.
\end{equation}
A solution of (\ref{eq:Xbint}) is any function $X(t)$ for which
both sides of (\ref{eq:Xbint}) exist and are equal.  Every solution
of (\ref{eq:Xb}) is a solution of (\ref{eq:Xbint}) but the converse
need not hold; whereas $X(t)$ must be differentiable for the left-hand
side of (\ref{eq:Xb}) to make sense, there is no explicit differentiability
requirement in (\ref{eq:Xbint}), only the requirement that
$b(s,X(s))$ be integrable.

The same idea carries over to random processes; although Brownian
motion cannot be differentiated, it can be used as an integrator,
and the state-space equations (\ref{eq:ctsse}) and (\ref{eq:ctoe})
can be written rigorously as functional equations by integrating
both sides.  However, there is in general no \textit{unique} way
of defining the integral of a process.

Since it may come as a surprise that different definitions of
integrals can give different answers, this phenomena will be
illustrated in the deterministic setting by considering how to
integrate H\"older continuous functions of order $\tfrac12$.  Such
functions are continuous but may oscillate somewhat wildly.

\noindent\textit{Technicality:}
Although sample paths of Brownian motion are almost surely not
H\"older continuous of order $\tfrac12$, they are almost surely H\"older
continuous of any order less than $\tfrac12$.  This is irrelevant here
because, as explained below, the relevant fact about Brownian motion
is that $\E{ |B(t+\delta t) - B(t)| }$ is proportional to $\sqrt{\delta
t}$ rather than $\delta t$.  The appearance of the expectation
operator characterises the stochastic approach to integration.
Interestingly, a complementary theory known as rough paths has been developed
recently~\cite{Lyons:1998bd, Lejay:2003df, bk:Friz:stochastic, Davie:2008js},
based partially on observations made in~\cite{Yamato:1979jz,
Krener:1980ih} and~\cite{Sussmann:1978wo}.

Recall the Riemann-Stieltjes integral $\int_0^1 f(s) \, dg(s)$
which, for smooth $f$ and $g$, is a limit 
of Riemann sums:
\begin{equation}
\label{eq:ifg}
\int_0^1 f(s)\,dg(s) = \lim_{N \rightarrow \infty}
\sum_{k=0}^{N-1} f(\tfrac{k}N) \left( g(\tfrac{k+1}N) - g(\tfrac{k}N) \right).
\end{equation}
The right-hand side of (\ref{eq:ifg}) gives the same answer 
if the right, not left, endpoint is used for each interval
$[\frac{k}N,\frac{k+1}N]$.  Indeed, the difference between using
left or right endpoints is
\begin{equation}
e_N = 
\sum_{k=0}^{N-1} \left( f(\tfrac{k+1}N) - f(\tfrac{k}N) \right)
    \left( g(\tfrac{k+1}N) - g(\tfrac{k}N) \right).
\end{equation}
If $f$ is smooth then $| f(\tfrac{k+1}N) - f(\tfrac{k}N) | \leq \alpha_f
\tfrac1N$ for some constant $\alpha_f \in \reals$, and
analogously for $g$.  Therefore, $| e_N | \leq \alpha_f \alpha_g
\sum_{k=0}^{N-1} \tfrac1N\,\tfrac1N$ and converges to zero as
$N \rightarrow \infty$.

If now $f$ and $g$ are not differentiable, but merely H\"older
continuous of order $\tfrac12$, then 
$| f(\tfrac{k+1}N) - f(\tfrac{k}N) | \leq \alpha_f | \tfrac1N
|^{\tfrac12}$ for some constant $\alpha_f \in \reals$, and analogously
for $g$, leading to $| e_N | \leq \alpha_f \alpha_g \sum_{k=0}^{N-1}
|\tfrac1N|^{\tfrac12}\,|\tfrac1N|^{\tfrac12}$ which converges to
$\alpha_f \alpha_g$, and not zero, as $N \rightarrow \infty$.  This
means it is possible for two different values of the integral
$\int_0^1 f(s)\,dg(s)$ to be obtained depending on whether
$f(\tfrac{k}N)$ or $f(\tfrac{k+1}N)$ is used in (\ref{eq:ifg}).

If at least one of $f$ or $g$ is smooth and the other is H\"older
continuous of order $\tfrac12$ then once again $|e_N| \rightarrow 0$.

If $g$ is replaced by real-valued Brownian motion $B(t)$ and the
above calculations carried out, a relevant quantity is the rate at
which the expected value of $|B(t+\delta t) - B(t)|$ decays to zero.
Since $B(t+\delta t) - B(t) \sim N(0,\delta t)$, $\E{ |B(t+\delta
t) - B(t)| }$ is proportional to $\sqrt{\delta t}$, analogous to
$|f(t+\delta t) - f(t)| \propto \sqrt{\delta t}$ for H\"older
continuous functions of order $\tfrac12$.  Not surprisingly then,
differences can appear for integrals of the form
$\int X(t)\,dY(t)$ when $X(t)$ and $Y(t)$ are 
stochastic processes; other integrals, such as $\int X(t)\,dt$ and $\int
h(t)\,dY(t)$, with $h$ smooth, are unambiguous. (This presupposes
$X(t)$ and $Y(t)$ are semimartingales~\cite{Protter:2004wf}.)

\noindent\textit{Technicalities:}
Lebesgue-Stieltjes theory requires finite variation, excluding
Brownian motion as an integrator.  H\"older continuous functions
of order greater than $\tfrac12$ can be integrated with respect to Lipschitz
functions using the Young integral without needing finite
variation~\cite{Young:1936ga}. Brownian motion falls just outside
this condition. Stochastic integration theory depends crucially on
integrators having nice \emph{statistical} properties for Riemann-sum
approximations to converge; see Section~\ref{sec:nacs}.  The sums
are sensitive to second-order information (cf.\@ \ito's
formula~\cite{bk:Oksendal:sde}), hence ``second-order calculus''
is fundamental to stochastic geometry~\cite[Section VI]{Emery:1989cy}.

\subsection{\ito\ and Stratonovich Integrals: An Overview}
\label{sec:oview}

Non-equivalent definitions of stochastic integrals~\cite{Protter:1979vz,
Russo:1993io} all involve taking limits of approximations but differ
in the approximations used and the processes that are allowed to
be integrators and integrands.  The two most common stochastic
integrals are the \ito\ and Stratonovich integrals.

The \ito\ integral leads to a rich probabilistic theory based
on a class of processes known as semimartingales, and a resulting
stochastic analysis that, in some ways, parallels functional analysis.
A tenet of analysis is that properties of a function $f(t)$
can be inferred from its derivative; for example,
a bound on $f(T)$ can be derived from bounds on $f'(t)$ because
$f(T) = \int_0^T f'(t)\,dt$.  (It is remarkable how often it
is easier to study an \emph{infinite} number of \emph{linear}
problems, namely, examining $f'(t)$ for each and every $t$ in the
range $0$ to $T$.)
Thinking of the random variable $X(T)$ as an infinite sum
of its infinitesimal differences --- $X(T) = \int_0^T dX(t)$ ---
suggests that by understanding the limiting behaviour of $X(t+\delta
t)-X(t)$ as $\delta t \rightarrow 0$, it is possible to infer
properties of $X(T)$ that may otherwise be difficult to infer
directly.

\noindent\textit{Technicality:}
If the decomposition $Y(T) = \int_0^T dY(t)$ of $Y(t) = f(X(t))$
is sought, the \ito\ formula allows the limiting behaviour of
$Y(t+\delta t) - Y(t)$ to be determined directly from the limiting
behaviour of $X(t+\delta t) - X(t)$.

The \ito\ integral does not respect geometry~\cite{Emery:1989cy};
it does not transform ``correctly'' 
to allow a coordinate-independent definition.
Nor does the \ito\
integral respect polygonal approximations: if $Y^k(t)$ is a sequence
of piecewise linear approximations converging to $Y(t)$ then it is
\emph{not} necessarily true that $\int X(t)\,dY^k(t) \rightarrow
\int X(t)\,dY(t)$.

The Stratonovich integral lacks features of the \ito\ integral but
respects geometry and polygonal approximations, making it suitable
for stochastic geometry and modelling physical phenomena.  Fortunately,
it is usually possible to convert from one to the other by adding
a correction term, affording freedom to choose the simpler for the
calculation at hand.

By respecting geometry, the Stratonovich integral behaves like its
deterministic counterpart; this is the
\textit{transfer principle}.  The archetypal example 
is that the trajectory $X(t)$ of a Stratonovich stochastic
differential equation $dX(t) = f(t,X(t)) \circ dB(t)$ stays on
a manifold $M$ if $f(t,X(t))$ lies in the tangent space $T_{X(t)}M$;
cf.\@ (\ref{eq:Xb}).

In terms of modelling, the \ito\ integral is suitable for approximating
inherently discrete-time systems by continuous-time systems (e.g.,
share trading), while the Stratonovich integral is suited to
continuous-time physical processes because it describes the limit
of piecewise smooth approximations~\cite{Mazo:2002vm}.

\noindent\textit{Technicality:}
On a non-Riemannian manifold, a Stratonovich integral but not an
\ito\ integral can be defined because only the former respects
geometry.  Only once a manifold is endowed with a connection can
an \ito\ integral be defined.

\subsection{Semimartingales and Adapted Processes}
The class of processes called semimartingales emerged over time by
attempts to push \ito's theory of integration to its natural limits.
Originally defined as ``adapted c\`adl\`ag processes decomposable
into the sum of a local martingale and a process of finite variation'',
the Bichteler-Dellacherie theorem \cite{Bichteler:1979jr,
Beiglboeck:2012wv} states that semimartingales can be defined
alternatively (in a simple and precise way~\cite{Protter:1986fj})
as the largest class of ``reasonable integrators'' around which can
be based a powerful and self-contained theory of \ito\ stochastic
integration~\cite{Protter:2004wf}.

\noindent\textit{Note:}
C\`adl\`ag and c\`agl\`ad processes~\cite{Protter:2004wf}
generalise continuous processes by permitting
``well-behaved'' jumps.

\noindent\textit{Technicality:}
By restricting the class of integrands, the class of integrators
can be expanded beyond semimartingales, leading to an integration
theory for fractional Brownian motion, for example.  Nevertheless,
the \ito\ theory remains the richest.

From an engineering perspective, primary facts are: all L\'evy
processes~\cite{Applebaum:2004to}, including the Poisson process
and Brownian motion, are semimartingales, as are all (adapted)
processes with continuously differentiable sample paths, and if a
semimartingale is passed into a system modelled by an \ito\ integral,
the output will also be a semimartingale.  From an analysis
perspective, semimartingales are analogues of differentiable functions
in that a process $X(t)$ can be written, and studied, in terms of
its differentials $dX(t)$; see Section~\ref{sec:oview}.  (Whereas
the differential of a smooth function captures only first-order
information, the Schwartz principle~\cite[(6.21)]{Emery:1989cy} is
that $dX(t)$ captures both first-order and second-order information.
This is \emph{stochastic} information though; sample paths of
semimartingales need not be differentiable.)

Crucial to \ito's development of his integral was the restriction
to \emph{adapted} processes: in $\int X\,dY$, \ito\ required $X(\tau)$
not to depend on $Y(t)$ for $t \geq \tau$.
(The borderline case $t = \tau$ results in \ito's requirement that
$X$ be c\`agl\`ad and $Y$ c\`adl\`ag.) A filtration formalises what
information has been revealed up to any given point in time.
Adaptedness to the filtration at hand is a straightforward technical
condition~\cite{bk:Oksendal:sde} taken for granted, and hence largely
ignored, in this primer.

\subsection{The \ito\ and Stratonovich Integrals}
\label{sec:itoint}

Let $X(t)$ be a continuous process and $Y(t)$ a semimartingale
(both adapted to the same filtration).  The \ito\ integral
\begin{equation}
\label{eq:origito}
Z(t) = \int_0^t X(s)\,dY(s)
\end{equation}
can be interpreted as a system with transfer function $X(t)$ that
outputs the semimartingale $Z(t)$ in response to the input $Y(t)$.

By~\cite[Theorem II.21]{Protter:2004wf},
(\ref{eq:origito}) is the limit (in probability) 
\begin{equation}
\label{eq:itodef}
Z(t) = \lim_{N \rightarrow \infty} \sum_{k=0}^{N-1}
    X(t_k) \big( Y(t_{k+1})-Y(t_k) \big)
\end{equation}
where $t_k = \tfrac{k}Nt$.  (The \ito\ integral naturally extends
to adapted c\`agl\`ad integrands.  This suffices for stochastic
differential equations.  With effort, further extensions are
possible~\cite{Protter:2004wf}.)

The Stratonovich integral~\cite{Stratonovich:1966tn}
is denoted
\begin{equation}
Z(t) = \int_0^t X(s)\circ dY(s).
\end{equation}
It can be thought of as
the limit (in probability) 
\begin{equation}
\label{eq:Zstr}
Z(t) = \lim_{N \rightarrow \infty} \sum_{k=0}^{N-1}
    X(\tau_k) \big( Y(t_{k+1})-Y(t_k) \big)
\end{equation}
where $t_k = \tfrac{k}Nt$ and $\tau_k = \frac{t_k+t_{k+1}}2$.
Alternatively, it can be evaluated as a limit of ordinary integrals.
Define the piecewise linear approximation
\begin{equation}
Y^\delta(k\,\delta +t)
= (1-t)\, Y(k\,\delta) + t\, Y((k+1)\,\delta )
\end{equation}
for $t \in [0,1]$ and non-negative integers $k$.
(For notational convenience, $\delta$ is being used here instead of
$\delta t$.)
Then
\begin{equation}
\label{eq:Zdel}
Z^\delta(t) = \int_0^t X(s)\,dY^\delta(s)
\end{equation}
is a well-defined ordinary integral, and $Z^\delta(t) \rightarrow
Z(t)$ as $\delta \rightarrow 0$. 
By differentiating $Y^\delta(s)$, (\ref{eq:Zdel}) becomes
\begin{equation}
Z^\delta(t) = \sum_{k=0}^{N-1}
    \tfrac{Y((k+1)\delta) - Y(k\delta)}{\delta}
    \int_{k\delta}^{(k+1)\delta} X(s)\,ds
\end{equation}
where $\delta = \tfrac{t}N$.  This agrees in the limit with (\ref{eq:Zstr})
whenever $\frac1\delta \int_{k\delta}^{(k+1)\delta} X(s)\,ds
\rightarrow X(\tau_k)$ where $\tau_k = (k+\tfrac12)\delta$.

\noindent\textit{Note:}
In the one-dimensional case, other reasonable approximations can
be used.  However, in general (i.e., when non-commuting vector
fields are involved), failure to use piecewise linear approximations
can lead to different answers~\cite{Cohen:2000hp}.

\noindent\textit{Technicality:}
Unlike for the \ito\ integral, it is harder to pin down conditions
for the Stratonovich integral to exist.  This makes the definition
of the Stratonovich integral a moving target.  It is generally preferable
to use (\ref{eq:Zstr}) with $X(\tau_k)$ replaced by $\tfrac12(X(t_k) +
X(t_{k+1}))$; the modified sum apparently converges under milder
conditions~\cite{Yor:1977ud}.  Alternatively, \cite{Nualart:2010we}
declares the Stratonovich integral to exist if and only if the
polygonal approximations with respect to all (not necessarily
uniform) grids converge in probability to the same limit, the limit
then being taken to be the value of the integral.  For reasonable
processes though, these definitions coincide.

\subsection{Stochastic Differential Equations}

The stochastic differential equation
\begin{equation}
dX(t) = b(t,X(t))\,dt + \Sigma(t,X(t))\,dB(t)
\end{equation}
is shorthand notation for the functional equation
\begin{equation}
\label{eq:sde}
X(t) = X(0) + \int_0^t b(s,X(s))\,ds +
    \int_0^t \Sigma(s,X(s))\,dB(s)
\end{equation}
which asks for a semimartingale $X(t)$ such that both sides of
(\ref{eq:sde}) are equal.
In higher dimensions, the \emph{diffusion coefficient} $\Sigma$ is
a function returning a matrix and the \emph{drift} $b$ a function returning
a vector.

The \ito\ equation (\ref{eq:sde}) can be solved
numerically by a standard forward-Euler method (known in the
stochastic setting as the Euler-Maruyama method)~\cite{Rao:1974wa,
Rumelin:1982jz, Kloeden:1992vt, Platen:1999cj, Iacus:2008bu}.  This
validates the otherwise \textit{ad hoc} models developed in
Section~\ref{sec:ssm}.

Replacing the \ito\ integral by the Stratonovich integral results
in a Stratonovich stochastic differential equation.  Since the
summation in (\ref{eq:Zstr}) involves evaluating $X(t)$ at a point
in the future --- the midpoint of the interval rather than the start
of the interval --- solving a Stratonovich equation necessitates
either using an implicit integration scheme (such as a predictor-corrector
method) or converting the Stratonovich equation into an \ito\
equation by adding a correction to the drift
coefficient~\cite{Rumelin:1982jz, Newton:1991dj, Kloeden:1992vt,
Platen:1999cj, Iacus:2008bu}.  For example, in the one-dimensional
case, solutions of the Stratonovich equation
\begin{equation}
dX(t) = \sigma(t,X(t))\circ dB(t)
\end{equation}
correspond to solutions of the \ito\ equation
\begin{equation}
dX(t) = \tfrac12 \sigma(t,X(t))\,\frac{\partial\sigma}{\partial X}(t,X(t))\,dt
+ \sigma(t,X(t))\,dB(t).
\end{equation}

\subsection{\ito\ Diffusions}

A process $X(t)$ is an \ito\ diffusion if it can be written in the
form (\ref{eq:sde}) where $B(t)$ is Brownian motion.  (In fact,
\ito\ developed his stochastic integral in order to write down
directly continuous Markov processes based on their infinitesimal
generators~\cite{Kuo:2006uk}.) Stratonovich equations of the analogous
form are also \ito\ diffusions.  \ito\ diffusions are particularly
nice to work with as they have a rich and well-established theory.

\noindent\textit{Technicality:}
The class of continuous Markov processes should be amongst the
simplest continuous-time processes to study since local behaviour
presumably determines global behaviour.  Surprisingly then, even
restricting attention to strongly Markov processes (thus ensuring
the Markov property holds also at random stopping times) does not
exclude complications, especially in dimensions greater than
one~\cite[Section IV.5]{Ikeda:1989wt}.  The generic term ``diffusion''
is used with the aim
of restricting attention to an amenable subclass of continuous
strongly Markov processes.  Often this results in diffusions being
synonymous with \ito\ diffusions but sometimes diffusions are 
more general.

\section{From Euclidean Space To Manifolds}
\label{sec:sdem}

Section~\ref{sec:ssm} used piecewise geodesic approximations to
define state-space models on manifolds.  These approximations
converge to solutions of stochastic differential equations.
Conversely, one way to understand Stratonovich differential equations
on manifolds is as limits of ordinary differential equations on
manifolds applied to piecewise geodesic approximations of sample
paths~\cite[Theorem 7.24]{Emery:1989cy}.

Stratonovich equations on a manifold $M$ can be defined using only
stochastic calculus in Euclidean space.  A class of candidate
solutions is needed: an $M$-valued process $X_t$ is a semimartingale
if the ``projected'' process $f(X_t)$ is a real-valued semimartingale
for every smooth $f\colon M \rightarrow \reals$; a shift in focus
henceforth makes the notation $X_t$ preferable to $X(t)$.  A
Stratonovich equation on $M$ gives a rule for $dX_t$ in terms of
$X_t$ and a driving process.  Whether a semimartingale $X_t$ is a
solution can be tested by comparing $df(X_t)$ with what it should
equal according to the stochastic chain rule; cf.\@ (\ref{eq:euc}).
If this test passes for every smooth $f\colon M \rightarrow \reals$
then $X_t$ is deemed a solution~\cite[Definition 1.2.3]{Hsu:2002tz}.

The following subsections touch on other ways of working with
stochastic equations on manifolds.

\subsection{Local Coordinates}

Manifolds are formally defined via a collection of charts composing
an atlas~\cite{Munkres:1991vu}.  (Charts were not mentioned in
Section~\ref{sec:rm} but exist by the implicit function theorem.)
A direct analogy is the representation of the Earth's surface by
charts (i.e., maps) in a cartographer's atlas.  Working in local
coordinates means choosing a chart --- effectively, a projection
of a portion of the manifold onto Euclidean space --- and working
with the projected image in Euclidean space using Cartesian
coordinates.  This is the most fundamental way of working with
functions or processes on manifolds.

A clear and concise example of working with stochastic 
equations in local coordinates is~\cite{Solo:2009ku}.  When written
in local coordinates, stochastic differential equations on manifolds
are simply stochastic differential equations in Euclidean space.

Working in local coordinates requires ensuring consistency across
charts.  When drawing a flight path in a cartographer's atlas and
nearing the edge of one page, changing to another page necessitates
working with both pages at once, aligning them on their overlap.
For theoretical work, this is usually only a minor inconvenience.
The inconvenience may be greater for numerical work.

\subsection{Riemannian Exponential Map}

On a general manifold there is no distinguished choice of local
coordinates; all compatible local coordinates are equally valid.
On a Riemannian manifold though, the Riemannian exponential map
gives rise to distinguished local parametrisations of the manifold,
the inverses of which are local coordinate systems called normal
coordinates.  The Riemannian exponential map has a number of
attractive properties.  It plays a central role in the smooth
approximation approach to stochastic development, and has been
used throughout this primer to generalise Brownian motion and
stochastic differential equations to manifolds.

The main disadvantage of the Riemannian exponential map is that,
in general, it cannot be evaluated in closed form, and its numerical
evaluation may be significantly slower than working
instead with extrinsic coordinates or even other local coordinate
systems.

\subsection{Extrinsic Coordinates}

Stochastic equations on manifolds embedded in Euclidean space can
be written as stochastic equations in Euclidean space that just so
happen to have solutions lying always on the manifold; see
Section~\ref{sec:seom}.

Even though all manifolds are embeddable in Euclidean space, there
is no reason to expect an arbitrary manifold to have an embedding
in Euclidean space that is convenient to work with (or even possible
to describe).  The dimension of the embedding space may be considerably
higher than the dimension of the manifold, making for inefficient
numerical implementations.  Numerical implementations may be prone
to ``falling off the manifold'', or equivalently, incur increased
computational complexity by projecting the solution back to
the manifold at each step.  (Falling off the manifold is only a
numerical concern because the transfer principle makes it easy to
constrain solutions of Stratonovich equations in $\reals^n$ to lie
on a manifold $M \subset \reals^n$; see Section~\ref{sec:oview}.)

\noindent\textit{Technicality:}
Given the use of Grassmann manifolds in signal processing
applications, it is remarked there is an embedding of Grassmann
manifolds into matrix space given by using orthogonal projection
matrices to represent subspaces.

\subsection{Intrinsic Operations}

Differential geometry is based on coordinate independence.  On a
Riemannian manifold, the instantaneous velocity and acceleration
of a curve have intrinsic meanings understandable without reference
to any particular coordinate system.  Similarly, stochastic operators
can be defined intrinsically~\cite{Emery:1989cy}.

With experience, working at the level of intrinsic operations is
often the most convenient and genteel.  It is generally the least
suitable for numerical work.

\section{A Closer Look at Integration}
\label{sec:acl}

In probability theory, the set of possible outcomes is denoted
$\Omega$, and random variables are (measurable) functions from
$\Omega$ to $\reals$.  Once Tyche, Goddess of Chance, decides the
outcome $\omega \in \Omega$, the values of all random variables
lock into place~\cite{Williams:1991wk}.  Therefore, for clarity,
stochastic processes will be written sometimes as $X_t(\omega)$.

Several issues are faced when developing
a theory for stochastic integrals
\begin{equation}
\label{eq:int}
Z(\omega) = \int_0^1 X_t(\omega)\,dY_t(\omega).
\end{equation}
(This generalises to $Z_t = \int_0^t X_t\, dY_t$ giving a
process $Z_t$ instead of a random variable $Z$.) That $X_t$ must
be predictable is already explained clearly in the literature; see
the section ``Na\"ive Stochastic Integration Is Impossible''
in~\cite{Protter:2004wf}.  Perhaps not so clearly explained in the
literature is the repeated claim that since sample paths of Brownian
motion have infinite variation, a pathwise approach
is not possible; a pathwise approach fixes $\omega$ and treats
(\ref{eq:int}) as a deterministic integral of sample paths.  Apparently
contradicting this claim are publications on pathwise
approaches~\cite{Karandikar:1995ct}!  This is
examined below in depth.

\subsection{Non-absolutely Convergent Series}
\label{sec:nacs}

The definition of the Lebesgue-Stieltjes integral $\int f(t)\,dg(t)$
explicitly requires the integrator $g(t)$ to have finite variation.
This immediately disqualifies (\ref{eq:int}) as a Lebesgue-Stieltjes
integral whenever the sample paths of $Y_t$ have infinite variation.
This is not a superfluous technicality;
for integration theory to be useful, integrals and limits must
interact nicely. At the very least, the bounded convergence theorem
should hold:
\begin{equation}
\label{eq:bct}
f^{(k)}(t) \rightarrow f(t) \implies
\int f^{(k)}(t)\,dg(t) \rightarrow \int f(t)\,dg(t)
\end{equation}
whenever the $f^{(k)}(t)$ are uniformly bounded.

Let $g(t)$ be a step function with $g(0)=0$ and transitions $g(t_n^+)
- g(t_n^-) = (-1)^{n+1} \frac1n$ at times $0 < t_1 < t_2 < \cdots
< 1$. Then $g(1) = \lim_{N \rightarrow \infty} \sum_{n=1}^N (-1)^{n+1}
\frac1n = \ln 2$.  Although declaring $\int_0^1 dg$ to equal
$g(1)-g(0)=\ln 2$ may seem the obvious choice, it or any other
choice would lead to a failure of the bounded convergence theorem:
the $f^{(k)}(t)$ in (\ref{eq:bct}) can be chosen so that the $\int
f^{(k)}(t)\,dg(t)$ are partial sums of $\sum_n (-1)^{n+1} \frac1n$
with summands rearranged, and since $\sum_n (-1)^{n+1} \frac1n$ is
not absolutely convergent, it can be made to equal any real number.
To achieve the limit $5$, start with $1+\frac13+\frac15+\cdots$
until $5$ is exceeded, then add $-\frac12-\frac14-\cdots$ until the
partial sum drops below $5$, and so forth.  Other orderings can
cause the partial sums to oscillate, preventing convergence.

To summarise, if $g(t)$ has infinite variation then the bounded
convergence theorem will fail because ``order matters''.

Let $Y_t(\omega)$ be a random step function with $Y_0(\omega) = 0$
and each increment $Y_{t_n^+} - Y_{t_n^-} = \pm \frac1n$ having
independent and equal chance of being positive or negative.  Every
sample path of $Y_t$ has infinite variation.  Given $\omega$, the
signs of each term in the sum $\sum_n \pm \frac1n$ become known,
and based on this knowledge, sequences can be constructed that
invalidate the bounded convergence theorem.

In applications, the order is reversed: a uniformly bounded
sequence $X_t^{(k)}$ is chosen with limit $X_t$, and of relevance
is whether the bounded convergence theorem $\lim_k \int X_t^{(k)}\,dY_t
= \int X_t\,dY_t$ holds for most if not all $\omega$.  Although a
sequence can be chosen to cause trouble for a single outcome $\omega$,
perhaps it is not possible for a sequence to cause trouble
simultaneously for a significant portion of outcomes. If $X_t$ can
depend arbitrarily on $\omega$ then trouble is easily caused for
all $\omega$, hence \ito's requirement for $X_t$ to be predictable.
To avoid this distraction, assume the $X_t^{(k)}$ are deterministic.

For $\sum_n (-1)^{n+1} \frac1n$ to converge to $5$, or even to
oscillate, requires long runs of positive signs followed by long
runs of negative signs.  Since the signs are chosen randomly,
rearranging summands may not be as disruptive as in the deterministic
case.  In fact, rearranging summands has \emph{no effect} on the
value of the sum for \emph{almost all} outcomes $\omega$.

\noindent\textit{Technicality:}
Define the random variables $A_n = \pm \frac1n$, $S_N = \sum_{n=1}^N
A_n$ and $T_N = \sum_{n=1}^N A_{\rho(n)}$ where $\rho$ is a permutation
of the natural numbers.  Put simply, the $S_N$ are the partial sums
of $\sum_n \pm \frac1n$ in the obvious order, and the $T_N$ are the
partial sums in some other order specified by $\rho$.  The $S_N$
form an $\mathcal{L}^2$-martingale sequence that converges almost
surely and in $\mathcal{L}^2$ to a finite random variable $S_\infty$;
see~\cite[Chapter 12]{Williams:1991wk}.  Similarly, $T_N$ converges
to a finite random variable $T_\infty$.  With a bit more effort,
and appealing to the Backwards Convergence Theorem~\cite{Protter:2004wf},
it can be shown that $S_\infty = T_\infty$ almost surely.

Since rearranging summands (almost surely) does not affect the
\emph{random} sum $\sum_n \pm\frac1n$ --- and the same holds for
the weighted sums $\sum_n \pm\frac1n X^{(k)}(t_n)$ --- the stochastic integral
$\int_0^1 dY_t(\omega)$ can be defined to be the random variable
$Y_1(\omega) - Y_0(\omega)$ without fear of the bounded convergence
theorem failing (except on a set of measure zero).  Here,
$Y_1(\omega) - Y_0(\omega)$ comes from the particular ordering
$\lim_{N \rightarrow \infty} \sum_{n=1}^N Y_{t_n^+} - Y_{t_n^-}$.

That Lebesgue-Stieltjes theory cannot define $\int_0^1 dY_t(\omega)$
as $Y_1(\omega) - Y_0(\omega)$ does not make this pathwise definition
incorrect.  The Lebesgue-Stieltjes theory requires absolute convergence
of sums, meaning all arrangements have the same sum.  The theory
of stochastic integration weakens this to requiring any two
arrangements have the same sum almost surely (or converge in some
other probabilistic sense).  The difference arises because there are
an infinite number of arrangements of summands and an infinite
number of outcomes $\omega$.  If $\omega$ is chosen first then two
arrangements can be found whose sums disagree, yet if two arrangements
are chosen first then their sums will agree for almost all $\omega$.

\subsection{Brownian Motion as an Integrator}
\label{sec:bmi}

The discussion in Section~\ref{sec:nacs} remains insightful when
$Y_t$ is Brownian motion.  The basic hope is that
(\ref{eq:int}) can be evaluated by choosing a sequence of partitions
$\{0=t_0^{(k)} < t_1^{(k)} < \cdots < t_{N^{(k)}}^{(k)} = 1\}$ with
mesh size $\sup_{i,j} |t_i^{(k)} - t_j^{(k)}|$ going to zero as $k
\rightarrow \infty$, and evaluating
\begin{equation}
\label{eq:hope}
\lim_{k \rightarrow \infty}
\sum_{n=0}^{N^{(k)}-1} X_{t_n^{(k)}} \left( Y_{t_{n+1}^{(k)}} - Y_{t_n^{(k)}}
\right).
\end{equation}
This can be interpreted as approximating $X_t$ by a sequence of
step functions, computing the integral with respect to each of these
step functions, and taking limits.  In particular, if different
partitions lead to different limits, the dominated convergence theorem
cannot hold.

A standard manipulation~\cite[p.\@ 30]{bk:Oksendal:sde} shows
(\ref{eq:hope}) simplifies to $\frac12 B_1^2 - \frac12
\sum_{n=0}^{N^{(k)}-1} ( B_{t_{n+1}^{(k)}} - B_{t_n^{(k)}} )^2$ for
the particular example $\int_0^1 B_t\,dB_t$.  If $B_t$ is Brownian
motion then the situation is somewhat delicate because the sample
paths are only H\"older continuous of orders strictly less than
$\tfrac12$; section~\ref{sec:detan}.  This manifests itself in the
existence of a partition such that~\cite{Morters:2010vm}
\begin{equation}
\label{eq:limsup}
\limsup_{k \rightarrow \infty} \sum_{n=0}^{N^{(k)}-1} \left(
B_{t_{n+1}^{(k)}} - B_{t_n^{(k)}} \right)^2 = \infty
\end{equation}
\emph{almost surely}, whereas the limit \emph{in probability} is
\begin{equation}
\label{eq:BinP}
\lim_{k \rightarrow \infty}
\sum_{n=0}^{N^{(k)}-1} \left( B_{t_{n+1}^{(k)}} - B_{t_n^{(k)}}
\right)^2 = 1,
\end{equation}
with (\ref{eq:BinP}) holding for every partition.  Therefore, it
appears stochastic integrals cannot be defined almost surely, but
at best in probability.  This is why, in (\ref{eq:itodef}), the
convergence is in probability.  Section~\ref{sec:path} discusses
the difference.

It is claimed in~\cite{Karandikar:1995ct} that almost sure convergence
is achievable; how is this possible without breaking the dominated
convergence theorem?  Crucially, the claim is that stochastic
integrals can be \emph{evaluated} using a particular sequence of
partitions yielding an almost sure limit.  \emph{The dominated
convergence theorem still only holds in probability.}

If (\ref{eq:hope}) converges almost surely then it is called
a pathwise evaluation of the corresponding stochastic integral.

\noindent\textit{Technicality:}
There are two standard ways of ensuring (\ref{eq:BinP}) converges
almost surely~\cite[Theorem I.28]{Protter:2004wf}.  The first is
to shrink the size of the partition sufficiently fast to satisfy
the hypothesis of the Borel-Cantelli lemma.  The second is to use
a nested partition, meaning no points are moved or deleted as $k$
increases, only new points added.  The same basic idea holds for
(\ref{eq:hope}) except care is required if the integrand
contains jumps; see for example~\cite[Theorem 2]{Karandikar:1995ct}.

\subsection{Pathwise Solutions of Stochastic Differential Equations}

There are several meanings of pathwise solutions of stochastic
differential equations.  Some authors mean merely a strong solution.
(Strong and weak solutions are explained in~\cite{bk:Oksendal:sde}.)
Or it could mean a limiting procedure converging almost
surely~\cite{Karandikar:1981vr}, in accordance with the definition
of pathwise evaluation of integrals in Section~\ref{sec:bmi}.  In
the strongest sense, it means the solution depends continuously on
the driving semimartingale~\cite{Sussmann:1978wo}, also known as a
robust solution.

Robust solutions are relevant for filtering.  Roughly speaking,
integration by parts applied to the Kallianpur-Striebel formula
leads to a desirable version of the conditional
expectation~\cite{Clark:2005ka, Bain:2009wo}.

In general, robust solutions will not exist for stochastic differential
equations with vector-valued outputs.  The recent theory of rough
paths demonstrates that continuity can be restored by only requiring
the output be continuous with respect to not just the input but
also integrals of various higher-order products of the
input~\cite{Lejay:2005tf, Lyons:1998bd, Lejay:2003df, bk:Friz:stochastic,
Davie:2008js}.

As explained in Section~\ref{sec:path}, not being able to find a
robust or even a pathwise solution is often inconsequential in
practice.

\noindent\textit{Technicality:}
This primer encourages understanding stochastic equations by thinking
in terms of limits of piecewise linear approximations.  Piecewise
linear interpolation of sample paths using nested dyadic partitions
are also good sequences to use in the theory of rough paths~\cite[Section
3.3]{Coutin:2007ve}.

\subsection{Almost Sure Convergence and Convergence In Probability}
\label{sec:path}

Let $Z_k(\omega)$ be a sequence of real-valued random variables.
Fixing $\omega$ means $Z_k(\omega)$ is merely a sequence of real
numbers.  If this sequence converges for almost all $\omega$ then
$Z_k$ converges almost surely.  A weaker concept 
is convergence in probability.  If $Z$ is a random
variable and, for all $\epsilon > 0$, $\lim_{k \rightarrow \infty}
\operatorname{Pr}(|Z_k - Z| > \epsilon) = 0$ then $Z_k$ 
converges to $Z$ in probability.

If $Z_k$ converges to $Z$ in probability, but not almost surely,
and $\omega$ is fixed, then $Z(\omega)$ cannot necessarily be
determined just from the sequence $Z_k(\omega)$.  Computing a limit
in probability requires knowing the behaviour of the $Z_k(\omega)$
as $\omega$ changes, explaining why the term ``pathwise'' is used
to distinguish almost sure convergence from convergence in probability.

This may lead to the worry that a device cannot be built whose
output is a stochastic integral of its input, e.g., a filter.
Indeed, only a single sample path is available in the real world!
However, a device will never compute $\lim Z_k(\omega)$ (unless a
closed-form expression is known), but rather, $Z_K(\omega)$ will
be used in place of $Z(\omega)$, where $K$ is sufficiently large.
And $Z_K(\omega)$ \emph{can} be computed from a single sample path.

In practice, what often matters is the mean-square error
$\E{(Z_K-Z)^2}$.  Neither almost sure convergence or convergence
in probability implies convergence in mean-square.  Importantly
then, convergence in mean-square can be used when defining stochastic
integrals~\cite{Wong:1985ky}.  (Convergence in mean-square implies
convergence in probability. Convergence in probability only
implies convergence in mean-square if the family of random variables
$|Z_k|^2$ is uniformly integrable.)

\noindent\textit{Note:}
A good understanding is easily gained by recalling first the
representative example of a random variable converging in probability
but not almost surely~\cite{bk:Stoyanov:1997} and then studying how
a bad partition satisfying (\ref{eq:limsup}) is constructed
in~\cite[p.\@ 365]{Morters:2010vm}.

\noindent\textit{Technicality:}
Although functional analysis plays an important role in probability
theory, some subtleties are involved because convergence in probability
is not a locally convex topology.  Also, there is no topology
corresponding to almost sure convergence~\cite{Ordman:1966cg}.

\subsection{Integrals as Linear Operators}

\newcommand{\Lp}{\mathbb{L}}
\newcommand{\Lv}{\mathcal{L}}

Ordinary integrals are nothing more than (positive) linear functionals
on a function space~\cite{Eberlein:1957wa}.  The situation is no
more complicated in the stochastic setting.

A stochastic integral $I\colon \Lp \rightarrow \Lv$ is a linear
operator sending a process $X_t$ to a random variable $Z = I(X_t)$,
more commonly denoted $Z =\int_0^1 X_t\,dB_t$ for some fixed process
$B_t$.  Here, $\Lp$ is a vector space of real-valued random processes
and $\Lv$ a vector space of real-valued random variables.

To be useful, $I$ should satisfy some kind of a bounded convergence
theorem.  This generally necessitates, among other things, restricting
the size of $\Lp$.  (In \ito's theory, $\Lp$ is usually the class
of adapted processes with c\`agl\`ad paths.)

Given suitable topologies on $\Lp$ and $\Lv$, one way to define $I$
is to define it first on a dense subset of $\Lp$ then extend by
continuity: if $X_t^{(k)} \rightarrow X_t$ then declare $I(X_t) = \lim_k
X_t^{(k)}$.  (Protter takes this approach, albeit when
$I$ maps processes to processes~\cite{Protter:1986fj}.)
This only works if the topology on $\Lp$ is sufficiently strong
that two different sequences converging to $X_t$ assign the same value
to $I(X_t)$.  Alternatively, restrictions can be placed on the allowable
approximating sequences; the value of $I(X_t)$ is evaluated as the
limit of $I(X_t^{(k)})$ but where there are specific guidelines on
the construction of the sequence $X_t^{(k)} \rightarrow X_t$.  In
both cases, the topology on $\Lv$ is also important; too strong and
$I(X_t^{(k)})$ need not converge.

If convergence in probability is used on $\Lv$ then there is the
perennial caveat that the limit is only unique modulo sets of measure
zero.  Although this means $I$ is not defined uniquely, any given
version of $I$ will still be a map taking, for fixed $\omega$, the
sample path $t \mapsto X_t(\omega)$ to the number $Z(\omega) =
I(X_t)(\omega)$.  In this respect, using the term pathwise to exclude
convergence in probability is misleading.  While a version of $I$
may not be constructable one path at a time in isolation, once a
version of $I$ is given, it can be applied to a single path in
isolation.  (The same phenomenon occurs for conditional expectation.)

\section{Stratonovich Equations on Manifolds}
\label{sec:seom}

Starting from this section, the level of mathematical detail
increases, leading to the study, in Section~\ref{sec:pe}, of an
estimation problem on compact Lie groups.

Let $M \subset \reals^n$ be a $d$-dimensional manifold.  Denote by
$B_t$ Brownian motion in $\reals^d$ with $B_t^i$ its $i$th element.
A class of processes on $M$ can be generated by injecting
Brownian motion scaled by a function $h$.  Precisely, the Stratonovich
stochastic differential equation
\begin{equation}
\label{eq:str}
X_t = \int_0^t \sum_{i=1}^d h_i(t,X_t) \circ dB_t^i
\end{equation}
defines a process $X_t$ that, if started on $M$, will remain on $M$
if the $h_i(t,X_t)$ always lie in $T_{X_t}M$.  In such cases, for any
smooth function $f\colon M \rightarrow \reals$, the projected process
$f(X_t)$ on $\reals$ satisfies
\begin{equation}
\label{eq:euc}
f(X_t) = \int_0^t \sum_{i=1}^d df\langle h_i(t,X_t) \rangle \circ dB_t^i
\end{equation}
where $df\langle v_p \rangle$ is the directional derivative of $f$
at $p$ in the direction $v_p \in T_pM$.  (If $f$ is locally the
restriction of a function $F\colon U \subset \reals^n \rightarrow \reals$
defined on an open subset of $\reals^n$ containing $p$ then $df\langle
v_p \rangle$ is the usual directional derivative of $F$ at
$p$ in the direction $v_p$.) 

\noindent\textit{Technicality:}
If $M$ was not embedded in $\reals^n$, so that (\ref{eq:str}) was
not a stochastic equation in Euclidean space, then (\ref{eq:euc})
could be used to define (\ref{eq:str}): $X_t$ is the semimartingale
on $M$ for which (\ref{eq:euc}) holds for every smooth function
$f$.

\section{Compact Lie Groups}
\label{sec:clg}

Lie groups~\cite{Stillwell:2008dj, Boothby:1986wz} are simultaneously
a manifold and a group. The group operations of multiplication and
inversion are smooth functions, thereby linking the two structures.
Extending an algorithm from Euclidean space to manifolds is often
facilitated by first extending it to compact Lie groups where the
extra structure helps guide the extension.

The simplest example of a compact Lie group is the circle $S^1$
with \emph{group multiplication} $(\cos\theta,\sin\theta) \cdot
(\cos\phi,\sin\phi) = (\cos(\theta+\phi),\sin(\theta+\phi))$ and
\emph{identity element} $(1,0) \in S^1 \subset \reals^2$.  Being a
subset of Euclidean space, the circle is seen to be compact by
being closed and bounded (Heine-Borel theorem).

\noindent\textit{Technicality:}
Compactness is a topological property.  In a sense, compact sets
are the next simplest sets to work with beyond sets containing only
a finite number of elements.  Compact sets are sequentially compact:
any sequence in a compact set has a convergent subsequence.

A more interesting example of a compact Lie group is the special
orthogonal group $SO(n)$ defined as the set of all $n \times n$
real-valued orthogonal matrices with determinant equal to one, with
matrix multiplication as group multiplication.  This implies
the identity matrix is the identity element.  By treating the
space of $n \times n$ real-valued matrices as Euclidean space ---
the Euclidean inner product is\footnote{The $\operatorname{vec}$
operator identifies the space $\reals^{n \times n}$
with the space $\reals^{n^2}$ by stacking the columns of a matrix
on top of each other to form a vector. Under this identification,
the Euclidean inner product on $\reals^{n^2}$ becomes the inner
product $\langle A,B\rangle = (\operatorname{vec} B)^T (\operatorname{vec}
A) = \Tr\{B^TA\}$.}
$\langle A, B \rangle = \Tr\{B^TA\}$ where $\Tr$ denotes the trace of
a matrix and superscript $T$ is matrix transpose --- it
can be asked if $SO(n)$ meets the requirements in Section~\ref{sec:rm}
of being a manifold (of dimension $\tfrac{n}2(n-1)$), to which the
answer is in the affirmative.  Furthermore, by considering elementwise
operations, it can be seen that \emph{group} multiplication and 
inversion, which are just \emph{matrix} multiplication and inversion, are
smooth operations.

The group $SO(3)$ has been studied extensively in the physics
and engineering literature~\cite{chirikjian:00}.  Studying $SO(3)$
is the same as studying rotations of three-dimensional space.

\subsection{The Matrix Exponential Map}
\label{sec:mem}

Choose a matrix $X \in \reals^{n \times n}$ satisfying $X^T X = I$
and $\det X = 1$.  In other words, choose an element $X$ of $SO(n)$.
Since $SO(n)$ is a group, it is closed under multiplication, hence
the sequence $X^0, X^1, X^2, X^3, \cdots$ lies in $SO(n)$, forming
a trajectory.  The closer $X$ is to the identity matrix, the closer
successive points in the trajectory become.  Note that $X$ is a
rotation of $\reals^n$, hence $X^n$ is nothing more than the rotation
$X$ applied $n$ times.

The idea of interpolating such a trajectory leads to the curve
$\gamma(t) = \exp(tA)$ where $\exp$ is the matrix
exponential~\cite{Moler:2003fn}.  Indeed, if $A \in \reals^{n \times
n}$ is such that $\exp(A)=X$ then $X^2 = \exp(2A)$, $X^3 = \exp(3A)$,
and so forth.  The set of $A$ for which $\exp(A)$ is an element of
$SO(n)$ forms a vector space; this nontrivial fact relates to $SO(n)$
being a Lie group.  Let $so(n)$ denote this set; it is the
Lie algebra of $SO(n)$.

Since $\exp(A)$ will lie in $SO(n)$ if $\exp(tA)$ lies in $SO(n)$
for an arbitrarily small $t > 0$, it suffices to determine $so(n)$
by examining the linearisation of $\exp$.  The constraint
$\exp(tA)^T\exp(tA)=I$ implies $(I+tA+\cdots)^T(I+tA+\cdots)=I$
from which it follows that $so(n) = \{A \in \reals^{n \times n}
\mid A+A^T=0\}$.  This shows also that $so(n)$ is the tangent space
of $SO(n)$ at the identity: $so(n) = T_ISO(n)$.

\subsection{Geodesics on $SO(n)$}
\label{sec:geod}

Recall from earlier that the Euclidean inner product on $\reals^{n
\times n}$ is $\langle A,B\rangle = \Tr\{B^TA\}$.  A geodesic on
$SO(n)$ is a curve with zero acceleration, where acceleration is
measured by calculating the usual acceleration in $\reals^{n \times
n}$ then projecting orthogonally onto the tangent space.

It may be guessed that $X^0, X^1, X^2, \cdots$ in Section~\ref{sec:mem}
lies on a geodesic, consequentially, $\gamma(t) = \exp(tA)$ should
be a geodesic whenever $A \in so(n)$.  The acceleration in Euclidean
space is $\gamma''(t) = \exp(tA)A^2$.  If this is orthogonal to the
tangent space of $SO(n)$ at $\gamma(t)$ then $\gamma$ is a geodesic.

Fix $t$ and let $Z = \gamma(t) = \exp(tA)$.  The tangent space
at $Z \in SO(n)$ is the vector space $Z\,so(n)$, that is, matrices
of the form $ZC$ with $C \in so(n)$.  This shows one of the many
conveniences of working with manifolds that are also groups; group
multiplication moves the tangent space at the identity element to
any other point on the manifold.  Then $\langle \gamma''(t), ZC
\rangle = \langle ZA^2,ZC \rangle = \langle A^2, C \rangle = 0$
whenever $A$ and $C$ are skew-symmetric matrices, proving $\gamma$
is indeed a geodesic.

A geodesic is completely determined by its initial position and
velocity, hence every geodesic starting at the identity element is
of the form $t \mapsto \exp(tA)$ for some $A \in so(n)$.

All geodesics starting at $X \in SO(n)$ are of the
form $t \mapsto X \exp(tA)$, that is, left multiplication by $X$
sends a geodesic to a geodesic.  Right multiplication also sends
geodesics to geodesics.  (This implies that for any $B \in so(n)$
and $X \in SO(n)$, there must exist an $A \in so(n)$ such that
$X\exp(tA) = \exp(tB)X$; indeed, $A = X^TBX$.)

If $\gamma(t) = X\exp(tA)$ then $\gamma'(0) = XA$.  Therefore,
the Riemannian exponential map is given by
\begin{equation}
\Exp_X(XA) = X\exp(A),\ X \in SO(n),\ A \in so(n).
\end{equation}

\subsection{Why the Euclidean Norm?}
\label{sec:eucnorm}

The matrix exponential map relates only to the group structure; $t
\mapsto \exp(tA)$ generates a one-parameter subgroup.  The Riemannian
exponential map relates only to the Riemannian geometry.  For there
to be a relationship between $\exp$ and $\Exp$ requires a careful
choice of Riemannian metric.

The special relationship between $SO(n) \subset \reals^{n \times
n}$ and the Euclidean inner product on $\reals^{n \times n}$ is
that the inner product is bi-invariant, meaning $\langle A, B \rangle
= \langle XA, XB \rangle = \langle AX,BX \rangle$ for any $A,B \in
so(n)$ and $X \in SO(n)$.  In words, the inner product between two
velocity vectors at the identity element of the Lie group is equal
to the inner product between the same two velocity vectors should
they be ``shifted'' to another point on the Lie group, either by
left multiplication or by right multiplication.  This is why both
left and right multiplication preserve geodesics. (Note: Since a
tangent vector is the velocity vector of a curve, the terms tangent
vector and velocity vector are used interchangeably.)

\noindent\textit{Technicalities:}
On a \emph{compact} Lie group, a bi-invariant metric can always be
found.  Different bi-invariant metrics lead to the same geodesics.
The most familiar setting to understand this in is Euclidean space;
changing the inner product alters the distance between two points,
but the shortest curve (and hence a geodesic) is still the same
straight line.  Similarly, on a product of circles, which is a
compact Lie group, a different scaling factor can be applied to
each circle to obtain a different bi-invariant metric, but the
geodesics remain the same.

\subsection{Coloured Brownian Motion on $SO(n)$} \label{sec:W}
\label{sec:cbm}

On an arbitrary Riemannian manifold, there is essentially only one
kind of Brownian motion.  On Lie groups (and other parallelisable
manifolds), the extra structure permits defining coloured
Brownian motion having larger fluctuations in some directions than
others.

If $G_t$ is a process on $SO(n)$,
it has both left and right
increments, namely, $G_{s+t}G_s^{-1}$ and $G_s^{-1}G_{s+t}$
respectively~\cite{Liao:2004dd}.
Non-commutativity implies these increments need not be
equal.  This leads to distinguishing left-invariant
Brownian motion from right-invariant Brownian motion.  As one is
nevertheless a mirror image of the other --- if $G_t$
is left-invariant Brownian motion then $G_t^{-1}$ is right-invariant
Brownian motion --- it generally suffices to focus on just one.

A process $G_t$ on $SO(n)$ is a left-invariant Brownian motion if
it has continuous sample paths and right increments that are
independent of the past and stationary.  (See~\cite[V.35.2]{Rogers:2000ur}
for details. Using \emph{right} increments leads to \emph{left}-invariance
of the corresponding transition semigroup~\cite{Liao:2004dd}.)
If it were stochastically anti-developed then a coloured
Brownian motion plus possibly a drift term would result.  

\noindent\textit{Technicality:}
In (\ref{eq:Wd}), the $\beta_{k,i}$ have zero mean, therefore, the
limiting process (\ref{eq:fbm}) has no drift.  This is equivalent
to restricting attention to left-invariant Brownian motions that
are also \emph{inverse invariant}~\cite[Section 4.3]{Liao:2004dd}.

An algorithm for simulating left-invariant Brownian motion on $SO(n)$
is now derived heuristically.  Unlike in Section~\ref{sec:two},
stochastic development is not used because left-invariant Brownian
motion relates to the group structure, not the geometric structure.
Nevertheless, the basic principle is the same: replace straight
lines by one-parameter subgroups.  The only issue is how to specify a
random velocity in a consistent way to achieve stationarity.
The velocity vector of $X\exp(tA)$ at $X$ is
$XA$.  Since a left-invariant process is sought, the velocity vector
$XA$ at $X$ should be thought of as equivalent to the velocity
vector $A$ at the identity (that is, left multiplication is used
to map $T_ISO(n)$ onto $T_XSO(n)$).  This leads to the following
algorithm.

Choose a basis $A_1,\cdots,A_d$ for $so(n)$ and let $\beta_{k,i}$
be a doubly-indexed sequence of \textit{iid} $N(0,1)$ Gaussian
random variables. Then a left-invariant Brownian motion starting
at $W_0 \in SO(n)$ is approximated by
\begin{equation}
\label{eq:Wd}
W_{(k+1) \delta} = W_{k\delta} \, \exp\left(\sqrt\delta\,\sum_{i=1}^d
    \beta_{k,i} A_i\right)
\end{equation}
for $k=0,1,2,\cdots$, where $\delta > 0$ is a small step size.
If necessary, these discrete points can be connected by geodesics to
form a continuous sample path.
That is, for an arbitrary $t$ lying between $k\delta$ and $(k+1)\delta$,
\begin{equation}
\label{eq:Wdt}
W_t = W_{k\delta}\,\exp\left(\frac{t-k\delta}{\sqrt\delta}
    \,\sum_{i=1}^d \beta_{k,i} A_i\right).
\end{equation}

The right increments $W_{k\delta}^{-1}W_{(k+1)\delta}$
are stationary and independent of the past values $\{W_t \mid t
\leq k\delta\}$, as required of left-invariant Brownian motion.

An outstanding issue is the effect of changing the basis $A_i$.
Comparing (\ref{eq:Wd}) with the corresponding formula for (white)
Brownian motion in Section~\ref{sec:two}, namely, $W_{(k+1)\delta}
= \Exp_{W_{k\delta}}\left( \sqrt\delta\,Z_k \right)$ where each
$Z_k$ is an ``$N(0,I)$'' Gaussian random variable on the tangent
space $T_{W_{k\delta}}SO(n)$, shows that (\ref{eq:Wd}) produces
white Brownian motion if the $A_i$ are an \emph{orthonormal} basis,
as to be expected.

Changing basis vectors in (\ref{eq:Wd}) is equivalent to changing
the colour of the driving Gaussian noise 
$\beta_{k,i}$.  Therefore, coloured Brownian motion with covariance
$C$ is defined as (the limit as $\delta \rightarrow 0$ of) the
process (\ref{eq:Wd})--(\ref{eq:Wdt}) where $\{A_i\}$ is a given
orthonormal basis and the \textit{iid} random vectors $\beta_k =
(\beta_{k,1}, \cdots,\beta_{k,d})$ have distribution $N(0,C)$.

\subsection{Formal Construction of Brownian Motion}
\label{sec:fbm}

Taking limits in (\ref{eq:Wd})--(\ref{eq:Wdt}) leads to a Stratonovich
equation for left-invariant Brownian motion on a Lie group $G$
of dimension $d$.  Let $A_1,\cdots,A_d$ be a basis for the Lie
algebra. Denote by $h_i$ the extension of $A_i$ to a left-invariant
vector field on $G$; for $G = SO(n)$, this is simply $h_i(X) =
XA_i$.  Let $B_t$ be Brownian motion on $\reals^d$.
Then the solution $W_t$ of
\begin{equation}
\label{eq:fbm}
W_t = \int_0^t \sum_{i=1}^d h_i(W_t) \circ dB_t^i
\end{equation}
is a left-invariant Brownian motion on $G$.

As in Section~\ref{sec:cbm}, if the $A_i$ are an \emph{orthonormal}
basis and the $B_t$ is coloured Brownian motion with $\E{B_tB_t^T}
= tC$ then $W_t$ is \emph{coloured Brownian motion with covariance $C$}.

\noindent\textit{Technicality:}
Whether $B_t$ is white Brownian motion and the $A_i$ changed, or
the $A_i$ are orthornormal and the colour of $B_t$ changed, is
a matter of choice.

\section{Brownian Distributions}
\label{sec:grv}

There is a two-way relationship between Gaussian random variables
and Brownian motion in $\reals$.  Brownian motion can be defined
in terms of Gaussian distributions.  Conversely, a Gaussian random
variable can be generated by sampling from Brownian motion $B_t$:
the distribution of $B_1$ is $N(0,1)$.

Lack of linearity prevents defining Gaussian random variables on
manifolds.  Nevertheless, sampling from Brownian motion on a
manifold produces a random variable that, if anything is, is the
counterpart of a Gaussian random variable. Such a random variable
is said to have a Brownian distribution.

Formally, given a Lie group $G$ of dimension $d$, an element $g \in
G$ and a positive semidefinite symmetric matrix $C \in \reals^{d
\times d}$, a random variable $X$ with (left) Brownian distribution
$N(g,C)$ has, by definition, the same distribution as $W_1$, where
$W_t$ is coloured Brownian motion with covariance $C$ started at
$W_0 = g$; see Section~\ref{sec:fbm}.  That is, coloured Brownian
motion started at $g$ is run for one unit of time and where it ends
up is a realisation of $N(g,C)$.

\section{Parameter Estimation}
\label{sec:pe}

Given a compact Lie group $G$, let $y_1,y_2,\cdots$ be \textit{iid} $N(g,C)$
distributed, as defined in Section~\ref{sec:grv}.  The aim is to
estimate $g$ and $C$ from the observations $y_i$.  In~\cite{Said:2012hi},
it was realised that by embedding $G$ in a matrix space $\reals^{n
\times n}$, an estimate of $g$ and (sometimes) $C$ is easily obtained from
the average of the images $Y_i \in \reals^{n \times n}$
of the $y_i \in G$ under the embedding.  This section re-derives
some of the results at a leisurely pace.  For full derivations and
consistency proofs, see~\cite{Said:2012hi}.

\subsection{Estimation on $SO(2)$}
\label{sec:sotwo}

Although $SO(2)$ is nothing other than the unit circle $S^1$ in
disguise, the theory below will be required later.

If $X$ is a $2 \times 2$ orthogonal matrix then its first column
corresponds to a point on the unit circle because the sum of the
squares of the two elements is one, and its second column, having
unit norm and being orthogonal to the first column, has only two
possibilities, corresponding to whether the determinant of $X$ is
$1$ or $-1$.  If $X \in SO(2)$, its second column is therefore fully
determined given its first column because the determinant of $X$
must be $1$.  (The relationship is deeper: $SO(2)$ is Lie group
isomorphic to $S^1$.)

Let $y \in SO(2)$ be randomly drawn from the distribution $N(g,\sigma^2)$
where the covariance matrix $C$ has been replaced by the scalar
$\sigma^2$ because the dimension of $SO(2)$ is $1$.  Since $SO(2)
\subset \reals^{2 \times 2}$, $y$ can be treated equally well as a
random matrix, and treating it as such, it is of interest to determine
its expected value, which in general will not lie in $SO(2)$.

Let $\D$ be the (infinitesimal) generator~\cite[Sec.\@ 7.3]{bk:Oksendal:sde}
of $W_t$, that is, for any smooth ``test function'' $f\colon G \rightarrow \reals$,
\begin{equation}
\label{eq:Df}
\D f(g) = \lim_{t \downarrow 0} \frac{\Eg{ f(W_t) } - f(g)}{t}
\end{equation}
where $W_t$ is Brownian motion starting at $g$ with covariance
matrix $C$; the superscript $g$ on the expectation symbol signifies
this dependence on $g$.  It follows from (\ref{eq:fbm}) and
(\ref{eq:euc}) in a well-known way that
\begin{equation}
\label{eq:DfB}
\D f = \frac12 \sum_{i,j=1}^d C_{ij} A_i^L A_j^L f
\end{equation}
where $A^L$ is the left-invariant vector field generated by $A$.
Vector fields act on functions by directional differentiation; $(A^L
f)(g) = \left.\frac{d}{dt}\right|_{t=0} f(g\exp(At))$.

If $G = SO(n)$ then $G$ is a subset of $\reals^{n \times n}$.  For
what follows, it suffices to assume $f$ is the restriction of a
linear function defined on $\reals^{n \times n}$.  In this case,
(\ref{eq:DfB}) simplifies to
\begin{equation}
\label{eq:lDfB}
\D f(g) = \frac12 \sum_{i,j=1}^d C_{ij} f(g A_i A_j).
\end{equation}
This will be derived from first principles for $SO(2)$. Let
\begin{equation}
A = \begin{bmatrix} 0 & -1 \\ 1 & 0 \end{bmatrix}
\end{equation}
be the chosen orthonormal basis for $so(2)$ with respect to the
scaled Euclidean inner product $\langle X, Y \rangle = \tfrac1{\sqrt2}
\Tr\{Y^TX\}$.  The scaling term is just a matter of convenience
(and convention).  Since $A$ is skew-symmetric and hence normal,
it is diagonalisable: $Q^{-1} A Q = D$ where
$D = \diag\{-\jmath,\jmath\}$ and
\begin{equation}
Q = \sqrt2 \begin{bmatrix} 1 & 1 \\ \jmath & -\jmath \end{bmatrix}.
\end{equation}
As $t \downarrow 0$,
(\ref{eq:Wd}) can be used to approximate $W_t$, yielding
\begin{equation}
\label{eq:aEg}
\Eg{W_t} \approx \int_{-\infty}^{\infty} g\, \exp(\sqrt{t} \beta A)\,
    p(\beta)\,d\beta
\end{equation}
where $p(\beta)$ is the probability density function of $N(0,\sigma^2)$.
By writing $\exp(\sqrt{t}\beta A) = Q\, \exp(\sqrt{t}\beta D)\, Q^{-1}$
and remembering the characteristic function of the distribution
$N(0,\sigma^2)$, it follows that
$\Eg{W_t} \approx \exp(-\tfrac{\sigma^2 t}2)\, g$.
Therefore,
\begin{equation}
\lim_{t \downarrow 0} \frac{\Eg{ W_t } - g}{t}
= (-\tfrac{\sigma^2}2)\,g.
\end{equation}
If $f$ is the restriction of a linear function on $\reals^{2 \times
2}$ then it can be interchanged with all the operations in
(\ref{eq:Df}), hence $\D f(g) = f( -\tfrac{\sigma^2}2\, g )
= -\tfrac{\sigma^2}2\,f(g)$.  Fortunately, this agrees with (\ref{eq:lDfB})
because $A^2=-I$, validating the approximation (\ref{eq:aEg}).

Time-homogeneity of (\ref{eq:fbm}) means $\D f(W_t)$ is
the instantaneous rate of
change of $\E{f(W_t)}$ for all $t$, not just $t=0$.  Roughly
speaking,
Kolmogorov's backward equation~\cite[Sec.\@ 8.1]{bk:Oksendal:sde}
evaluates $\E{f(W_t)}$ by integrating these infinitesimal changes.
Since the changes are path dependent, a PDE rather than an ODE is
required: $u(t,g) = \Eg{f(W_t)}$ satisfies
\begin{equation}
\label{eq:kbe}
\frac{\partial u}{\partial t} = \D u,\quad u(0,g) = f(g)
\end{equation}
where $\D$ acts on the function $g \mapsto u(t,g)$.
Although $f$ should map to $\reals$, 
letting it be the identity map on $\reals^{2 \times 2}$
does no harm, in which case the solution to (\ref{eq:kbe}),
now a PDE on $\reals^{2 \times 2}$, is
\begin{equation}
\Eg{W_t} = u(t,g) = \exp(-\tfrac{\sigma^2 t}2)\, g.
\end{equation}
It just so happens that this true solution agrees with the earlier
approximation for $\Eg{W_t}$.  In general this will not be the case
and solving the PDE will be necessary.

\noindent\textit{Technicality:}
Kolmogorov's backward equation is a PDE for evaluating $\tilde
u(t,g) = \E{f(W_T)\mid W_t=g}$ for a fixed $T$ by integrating
backwards from time $t=T$.  By time homogeneity, $\tilde u(t,g) =
\E{f(W_{T-t})\mid W_0=g} = u(T-t,g)$, explaining why (\ref{eq:kbe})
integrates forwards from time $T-T=0$.

In summary, it has been derived from first principles
that
\begin{equation}
\E{y} = \Eg{W_1} = \exp(-\tfrac{\sigma^2}2)\, g.
\end{equation}

Given $y_1,\cdots,y_m \in SO(2)$, the method of moments dictates
equating the sample mean with the true mean, namely, estimating $g$
and $\sigma^2$ by solving $\exp(-\tfrac{\sigma^2}2)\, g = y$ where $y =
\tfrac1m \sum_{i=1}^m y_i$ is the (extrinsic) sample average obtained
by treating the $y_i$ as $2 \times 2$ matrices.  As
$y$ will not necessarily be of the form $\exp(-\tfrac{\sigma^2}2)\, g$,
the equation should be solved in the least-squares sense, facilitated
by the polar decomposition of $y$.  If $y=UP$ where
$U$ is orthogonal and $P$ positive-definite symmetric then $g$
can be estimated by $U$ and $\exp(-\tfrac{\sigma^2}2)$ can be estimated
by $\tfrac12 \Tr\{P\}$.  Equivalently, $g$ can be estimated by the
``Q'' matrix in the QR decomposition of $y$.

\subsection{Estimation on $SO(3)$}

The same steps are followed as in Section~\ref{sec:sotwo}.
Choose the basis $A_1$, $A_2$, $A_3$ to be, respectively,
\begin{equation}
\label{eq:Athree}
\begin{bmatrix} 0&-1&0\\1&0&0\\0&0&0 \end{bmatrix},
\begin{bmatrix} 0&0&-1\\0&0&0\\1&0&0 \end{bmatrix},
\begin{bmatrix} 0&0&0\\0&0&-1\\0&1&0 \end{bmatrix}.
\end{equation}
If $f$ is linear on $\reals^{3 \times 3}$ then, from (\ref{eq:lDfB}),
\begin{equation}
\label{eq:Z}
\D f(g) = f(gZ),\quad
Z = \frac12\sum_{i,j} C_{ij} A_i A_j.
\end{equation}
As before, a solution to $\frac{\partial u}{\partial t} = \D u$
is $u(t,g) = g\exp(tZ)$; note $u$ is linear in $g$ for all $t$,
hence the validity of using (\ref{eq:Z}) which neglects
nonlinear terms of $f$.
In summary, if $y \sim N(g,C)$ then $\E{y} = g \exp(Z)$.

A polar decomposition of the sample average allows the estimation
of $g$ and $\exp(Z)$ because $\exp(Z)$ is a positive-definite
symmetric matrix, a consequence of $Z$ being symmetric (which in
turn is due to $C$ being symmetric and the $A_i$ anti-symmetric).
The only remaining question is if $C$ can be recovered from $Z$.
Direct expansion shows that for $SO(3)$,
\begin{equation}
Z = -\frac12 \begin{bmatrix}
C_{11}+C_{22}    & C_{23}  &  -C_{13} \\
C_{32}    & C_{11} +C_{33} & C_{12}  \\
-C_{31}    &   C_{21} &   C_{22} + C_{33}
\end{bmatrix}.
\end{equation}
Therefore, for $SO(3)$, the parameters $g$ and $C$ can be
estimated from a polar decomposition of the sample average.

\subsection{Estimation on $SO(n)$ for $n > 3$}
\label{sec:geq3}

The above ideas will not allow all elements of $C$ to be estimated
in higher dimensional cases for the simple reason that the dimension
of $SO(n)$ is $d = \tfrac{n}2(n-1)$, hence an $N(g,C)$ distribution on
$SO(n)$ requires $d + \tfrac{d}2(d+1) = \tfrac{n}8(n-1)(n^2-n+6)$ real-valued
parameters for its specification, which would exceed $n^2$, the
dimension of the ambient space $\reals^{n \times n} \supset SO(n)$,
if $n$ were larger than three.

Nevertheless, the calculations in the previous section remain valid
for $n > 3$.  The expected value of $y \sim N(g,C)$ is $g\exp(Z)$
where $Z$ is given by (\ref{eq:Z}).  In particular, $g$ can be
estimated as before, either via a polar decomposition or as the
``Q'' matrix in the QR decomposition of the sample average.
Furthermore, if $C$ is known to have a particular structure, the
simplest example being if $C$ is diagonal, then it may still be
possible to recover $C$ from $Z$.

\subsection{Estimation on a Compact Lie Group}
\label{sec:estimation}

The above ideas extend directly to $SU(n)$ and $U(n)$, the real
(not complex!) Lie groups consisting of complex-valued unitary
matrices, the former with the additional requirement that the
determinant be one.  See~\cite{Said:2012hi} for details.

Let $G$ be an arbitrary compact Lie group.  Group representation
theory studies (among other things) smooth maps of the form $f\colon G
\rightarrow U(n)$ which preserve group operations, meaning $f(g^{-1})
= [f(g)]^{-1}$ and $f(gh) = f(g)\,f(h)$.  Such maps are called
(unitary) representations.

Let $f$ be a representation of $G$ and $y \sim N(g,C)$, the latter
with respect to the orthonormal basis $A_1,\cdots,A_d$ for the Lie
algebra of $G$.  The same steps as before will be used to calculate
$\Eg{f(y)}$.  Although $f$ is not linear, it has another simplifying
feature instead; it preserves group operations.  In particular,
$A_i^L f(g) = \frac{d}{dt} f(ge^{tA_i}) = f(g) \frac{d}{dt} f(e^{tA_i})
= f(g) (A_if)$ where the intermediate derivatives are to be evaluated
at $t=0$ and $A_if$ denotes the derivative of $f$ at the identity
element in the direction $A_i$.  Note that $A_if$ will be an element
of the Lie algebra of $U(n)$, that is, a skew-Hermitian matrix.
Therefore,
\begin{equation}
\label{eq:Dfr}
\D f(g) = f(g) Z_f,\quad Z_f = \frac12\sum_{i,j} C_{ij}\,(A_if)\,(A_jf)
\end{equation}
where $Z_f$ is a symmetric matrix.  It is readily verified that
$u(t,g) = f(g)\exp(tZ_f)$ solves $\frac{\partial u}{\partial t} = \D u$.
(The map $g \mapsto u(t,g)$ preserves group multiplication which is
the only requirement for (\ref{eq:Dfr}) to be valid.) Thus,
\begin{equation}
\E{f(y)} = f(g) \exp(Z_f)
\end{equation}
where $Z_f$ is the matrix given in (\ref{eq:Dfr}).  Regardless of
the choice of $f$, $Z_f$ will be symmetric and so $\exp(Z_f)$ will
be a positive-definite symmetric matrix.

Once again, a polar decomposition can be used to estimate $f(g)$
and $\exp(Z_f)$.  How much of $g$ and $C$ can be recovered depends
on the mappings $g \mapsto f(g)$ and $C \mapsto Z_f$.

There is always (at least) one representation $f\colon G \rightarrow
U(n)$ that can be written down immediately for a compact Lie group
equipped with a bi-invariant metric, and that is the adjoint
representation.  Recall that the adjoint representation $\Ad\colon G
\rightarrow \Aut(\mathfrak{g})$ maps a group element $g \in G$ to
an automorphism of the corresponding Lie algebra $\mathfrak{g}$.
Given an orthonormal basis for $\mathfrak{g}$, any automorphism of
$\mathfrak{g}$ can be represented in matrix form with respect to
this basis.  A consequence of the basis being orthonormal and the
metric being bi-invariant is that these matrices will be unitary.
That is, by choosing an orthonormal basis, the adjoint representation
can be written in the form $f\colon G \rightarrow U(n)$ where $n$ is the
dimension of $G$.  This allows the application of the above extrinsic
averaging method to arbitrary compact Lie groups.

Furthermore, the adjoint representation (or any other representation)
can be used to augment any information already available about the
parameters $g$ and $C$.  That is to say, given \textit{iid} samples
$y_1,\cdots,y_m \sim N(g,C)$, in addition to forming $y = \tfrac1m\sum_i
y_i$ (assuming the $y_i$ belong to an appropriate Lie group such
as $SO(n)$), the average $\tfrac1m\sum_i f(y_i)$ can also be formed.
Its polar decomposition would yield an estimate of $f(g)$ and of
$Z_f$, thereby supplementing existing information about $g$ and
$C$.  This can be repeated for any number of representations
$f_1,f_2,\cdots$.

\section{A Few Words on Manifold-Valued Estimators}
\label{sec:afw}

Estimation theory extends to manifolds~\cite{Hendriks:1991hl,
Hendriks:1998bf, Pennec:2004ud}.  Several concepts normally taken
for granted, such as unbiasedness of an estimator, are not geometric
concepts and hence raise the question of their correct generalisations
to manifolds.

The answer is that the difficulty lies not with manifolds, but with
the absence of meaning to ask for an estimate of a parameter.  The
author believes firmly that asking for an estimate of a parameter
is, \textit{a priori}, a meaningless question.  It has been given
meaning by force of habit.  An estimate only becomes useful once
it is used to make a decision, serving as a proxy for the
unknown true parameter value.  Decisions include: the action taken
by a pilot in response to estimates from the flight computer; an
automated control action in response to feedback; and, what someone
decides they hear over a mobile phone (with the pertinent question
being whether the estimate produced by the phone of the transmitted
message is intelligible).  \emph{Without knowing the decision to
be made, whether an estimator is good or bad is unanswerable.} One
could hope for an estimator that works well for a large class of
decisions, and the author sees this as the context of estimation
theory.

The definition of unbiasedness in~\cite{Hendriks:1991hl} with respect
to a loss function accords with this principle of application
dependence.  Otherwise, two common but \textit{ad hoc} definitions
of the mean are the extrinsic mean and the Karcher mean.  The
extrinsic mean of a random variable $x \in M \subset \reals^n$
simply forgets the manifold and takes the mean of $x$ treated as
an $\reals^n$-valued random variable.  The extrinsic mean generally
will not lie on $M$.  The Karcher mean~\cite{Manton:2004vk,
Manton:2006vi} looks to be an ideal generalisation of the mean, but
on close inspection, problems arise if the manifold has positive
curvature and the samples are sufficiently far apart to prevent the
Karcher mean from existing or being unique.  (Since
manifolds are locally Euclidean, asymptotic properties of estimators
are more readily generalised to manifolds.)

In summary, the author advocates working backwards, starting from
a (representative) decision problem and deriving a suitable estimator
for making good decisions.  This viewpoint should clarify otherwise
arbitrary choices that must be made along the way.

\noindent\textit{Technicality:}
The James-Stein estimator helps demonstrate the definition of
a good estimator is application dependent~\cite{Manton:1998ha}:
remarkably, in terms of minimising the mean-square error, the best
estimate of $\mu \in \reals^3$ given $X \sim N(\mu,I)$ is not
$\widehat\mu = X$.  Since usually (but not always) the ``best''
estimator to use in practice really is $\widehat\mu = X$, this
illustrates a shortcoming of the default assumption that
mean-square error is \emph{always} a good measure of performance.

\section{Conclusion}
\label{sec:concl}
Pertinent background not easily accessible in the literature on how
to understand and work with stochastic processes on manifolds was
presented.  Stochastic development was introduced to ease the
transition from Euclidean space to manifolds: stochastic development
maps each sample path in Euclidean space to a sample path on a
manifold that preserves the infinitesimal properties of the process.
Not surprisingly then, Brownian motion on a manifold is the stochastic
development of Brownian motion in Euclidean space.  It was also
explained that stochastic development, Stratonovich integrals and
Stratonovich differential equations could all be understood as
limits of piecewise linear approximations in Euclidean space and,
more generally, as limits of piecewise geodesic approximations on
manifolds.

Lie groups are manifolds with a compatible group structure.  The
group structure permits defining ``coloured'' Brownian motion on
Lie groups.  Sampling from coloured Brownian motion produces random
variables with Brownian distributions that are generalisations of
Gaussian random variables.  Formulae from~\cite{Said:2012hi} were
re-derived for estimating the parameters of Brownian distributed
random variables.  The derivation demonstrated that standard
techniques, such as using Kolmogorov's backward equation, remain
applicable in the more general setting of processes on manifolds.

\section*{Acknowledgements}

Both the content and presentation of the paper have benefited from
earlier discussions with Dr Salem Said and Professor Victor Solo,
and from thorough critiques by reviewers.

\bibliographystyle{IEEEtran}
\bibliography{sdg}

%

\begin{IEEEbiographynophoto}{Jonathan Manton}
received the B.Sc.\@ degree in mathematics and the B.Eng.\@ degree
in electrical engineering in 1995 and the Ph.D.\@ degree in 1998,
all from the University of Melbourne.  He holds a Distinguished
Chair at the University of Melbourne with the title Future Generation
Professor. He is also an Adjunct Professor in the Mathematical
Sciences Institute at the Australian National University.
\end{IEEEbiographynophoto}

\end{document}